\documentclass[11pt,reqno]{amsart}
\usepackage{amssymb}
\usepackage{amsmath}
\usepackage{bm}
\usepackage{amsthm}
\usepackage{amsbsy}
\usepackage{psfrag}
\usepackage{graphics}
\usepackage{graphicx}
\usepackage{pgfplots}
\usepackage{relsize}
\usepackage{color}
\usepackage{xcolor}
\usepackage{pdfpages}
\usepackage{url}
\usepackage{hyperref}
\usepackage[square,sort&compress,comma,numbers]{natbib}
\usepackage{graphicx}
\usepackage{caption}
\usepackage{subcaption}
\theoremstyle{plain}
\newtheorem{theorem}{Theorem}[section]
\newtheorem{proposition}[theorem]{Proposition}

\newtheorem{lemma}[theorem]{Lemma}

\newtheorem{example}[theorem]{Example}
\newcommand{\norm}[1]{\left\Vert #1 \right\Vert}
\newtheorem{remark}[theorem]{Remark}

\begin{document}
\allowdisplaybreaks[4]
\numberwithin{figure}{section}

%
\numberwithin{table}{section}
 \numberwithin{equation}{section}
%

\title[FEM for Stokes Dirichlet Boundary Control  Problem]
{A two level finite element method for Stokes constrained Dirichlet boundary control problem}

\author[T. Gudi]{ Thirupathi Gudi}
 \address{Department of Mathematics, Indian Institute of Science, Bangalore - 560012, India}
 \email{gudi@iisc.ac.in}

 \author[R. C. Sau]{Ramesh Ch. Sau}

\address{Department of Mathematics, Indian Institute of Science, Bangalore - 560012, India}

\email{rameshsau@iisc.ac.in}

\date{}

\begin{abstract}
	In this paper we present a finite element analysis for a Dirichlet boundary
	control problem governed by the Stokes equation. The Dirichlet control is considered in a convex closed subset of the energy space $\mathbf{H}^1(\Omega).$
Most of the previous works on the Stokes Dirichlet boundary control problem deals with either tangential control or the case where the flux of the control is zero. This choice of the control is very particular and their choice
of the formulation leads to the control with limited regularity. To overcome this difficulty, we introduce the Stokes problem with outflow condition and the control acts on the Dirichlet boundary only hence our control is more general and it has both the tangential and normal components.
We prove well-posedness and discuss on the regularity of the control problem. The first-order optimality condition for the control leads to a Signorini problem. We develop a two-level finite element discretization by using $\mathbf{P}_1$ elements(on the fine mesh) for the velocity and the control variable and $P_0$ elements (on the coarse mesh) for the pressure variable. The standard energy error analysis gives $\frac{1}{2}+\frac{\delta}{2}$ order of convergence when the control is in $\mathbf{H}^{\frac{3}{2}+\delta}(\Omega)$ space. Here we have improved it to $\frac{1}{2}+\delta,$ which is optimal. Also, when the control lies in less regular space we derived optimal order of convergence up to the regularity. The theoretical results are corroborated by a variety of numerical tests.
\end{abstract}


\maketitle
\allowdisplaybreaks
\def\R{\mathbb{R}}
\def\I{\mathbb{I}}
\def\dx{{\rm~dx}}
\def\dy{{\rm~dy}}
\def\ds{{\rm~ds}}
\def\y{\mathbf{y}}
\def\u{\mathbf{u}}
\def\V{\mathbf{V}}
\def\w{\mathbf{w}}
\def\J{\mathbf{J}}
\def\f{\mathbf{f}}
\def\g{\mathbf{g}}
\def\div{{\rm div~}}
\def\v{\mathbf{v}}
\def\z{\mathbf{z}}
\def\b{\mathbf{b}}
\def\cV{\mathcal{V}}
\def\V{\mathbf{V}}
\def\cT{\mathcal{T}}
\def\sjump#1{[\hskip -1.5pt[#1]\hskip -1.5pt]}

\section{Introduction}\label{sec:Intro}
 In this paper, we consider the following Dirichlet boundary control problem governed by Stokes equation
 \begin{equation}\label{min:j}
 \min ~J(\mathbf{u},\mathbf{y})= \frac{1}{2}\norm{\mathbf{u}-\mathbf{u}_d}_{0,\Omega}^2+ \frac{\rho}{2}\norm{\nabla\mathbf{y}}_{0,\Omega}^2
 \end{equation}	
subject to,
\begin{equation}\label{stokes_eq_intro}
\begin{split}
-\Delta \mathbf{u}+\nabla{p}&=\mathbf{f} \quad \text{in}\;\Omega, \\
\nabla\cdot{\mathbf{u}}&=0 \quad \text{in}\; \Omega,\\
\bf{u}&=\mathbf{y} \quad \text{on}\; \Gamma_C,\\
\bf{u}&=\mathbf{0} \quad \text{on}\; \Gamma_D,\\
\frac{\partial \mathbf{u}}{\partial \mathbf{n}}-p \mathbf{n}&=\mathbf{0}\;\; \text{on} \;\; \Gamma_N,
\end{split}	
\end{equation}
with the control constraints  $$ \mathbf{y}_a\leq\bm{\gamma}_{0}(\mathbf{y})\leq \mathbf{y}_b\text{ on } \Gamma_C.$$

\noindent
In Section \ref{sec:ModelProblem}, we elaborate on the aforementioned problem. We propose a finite element approximation of the state and the control variable in order to discretize the above system. The first discussion on discretization of optimal control problems governed by partial differential equations was in the papers of Falk \cite{Falk:1973:Control}, Gevici \cite{Gevec:1979:Control}. Subsequently, many significant contributions have been made to this field. A control can act in the interior of a domain, in this case, we call distributed, or on the boundary of a domain, we call boundary (Neumann or Dirichlet) control problem. We refer to \cite{meyerrosch:2004:Optiaml,Ramesh:2018} for distributed control related problem, to \cite{CasasM:2008:Neum,Ramesh:2018} for the Neumann boundary control problem.

The Dirichlet boundary optimal control problems play an important
role in the context of the computational fluid dynamics, see, e.g., \cite{Fursikov,Hinze_Kunisch_2004}. The \textit{a priori} error analysis for such problems can be traced back to \cite{CR:2013:Dirich}.
The literature on Dirichlet boundary control problem outlines various approaches.
One typical method is to choose control from the $L^2(\partial\Omega)$-space:
 \begin{equation}\label{min:j01}
 \min_{y\in L^2(\Gamma)} ~J(u,y):= \frac{1}{2}\norm{u-u_d}_{0,\Omega}^2 + \frac{\rho}{2}\norm{y}_{0,\partial\Omega}^2,
 \end{equation}	
subject to Poisson problem
\begin{equation}\label{elliptic_pde}
\Delta u=f \;\;\text{in}\; \Omega,\quad
u=y \;\;\text{on} \;\;\partial\Omega.
\end{equation}
Due to the fact that the Dirichlet data is only in space $L^2(\partial\Omega)$, we need to understand the state equation \eqref{elliptic_pde} in an ultra-weak sense. This ultra weak formulation is easy to implement
and usually results in optimal controls with low regularity. Especially, when the problem is posed on a convex
polygonal domain, the control $y$ exhibits layer behaviour at the corners of the domain. This is because, it is determined
by the normal derivative of the adjoint state, whereas in a nonconvex polygonal domain the control may have singularity around a corner point, for more details one can see \cite{CR:2013:Dirich}. Another approach, as in \cite{CMR:2009:Dirich}, is the Robin boundary penalization which transforms the
Dirichlet control problem into a Robin boundary control problem.

One other popular approach is to find controls from the energy space, i.e., $H^{1/2}(\partial\Omega)$:
 \begin{equation}\label{min:j02}
 \min_{y\in H^{1/2}(\partial\Omega)} ~J(u,y):= \frac{1}{2}\norm{u-u_d}_{0,\Omega} + \frac{\rho}{2} |y|_{H^{1/2}(\partial\Omega)}^2,
 \end{equation}	
 we refer \cite{Steinbach:2014:Dirichlet} for this approach. We can define the standard weak solution for the state
 equation \eqref{elliptic_pde} with this choice of control. This approach introduces the Steklov-Poincar\'e operator to establish a new cost functional. The Steklov-Poincar\'e operator transforms the Dirichlet data into a Neumann data by using harmonic extension of the Dirichlet data;
  but this type of abstract operator may cause some difficulties in numerical implementation.
  It is well known that for a given $y\in H^{1/2}(\partial\Omega)$ there exists a harmonic extension $u_y\in H^1(\Omega)$ such that $|y|_{H^{1/2}(\partial\Omega)}$ can be equivalently defined as
  \begin{equation*}
  |y|_{H^{1/2}(\partial\Omega)}:=\norm{\nabla u_y}_{0,\Omega}.
  \end{equation*}
  This motivates to choose the control from $H^1(\Omega)$ space as in \eqref{min:j}, i.e.,

  \begin{equation}\label{min:j03}
  \min_{y\in H^1(\Omega)} ~J(u,y):= \frac{1}{2}\norm{u-u_d}_{0,\Omega}^2 + \frac{\rho}{2}\norm{\nabla y}_{0,\Omega}^2.
  \end{equation}

   This approach for Dirichlet boundary control problem was first introduced in the paper \cite{Gudi:2017:DiriControl}. Since the control is sought from the space $H^1(\Omega)$ they do not need the Steklov-Poincar\'e operator and hence this method is computationally very efficient.
     In the paper \cite{Gudi:2017:DiriControl}, one can find only unconstrained control, an improved analysis for constrained control can be found in \cite{Gudi_Ramesh:ESIAM_COCV}.
\par
 In this article, we consider a Stokes equation with mixed boundary conditions and the control acts on the Dirichlet boundary only.
 In the literature of Stokes Dirichlet control problem, we can see that two types of control are chosen. The first one is tangential control  i.e., the control acts only in the tangential direction of the boundary (see \cite{gong:HDG}). In \cite{gong:HDG} the authors propose hybridize discontinuous Galerkin (HDG) method to approximate the solution of a tangential Dirichlet boundary control problem with an $L^2$
 penalty on the boundary control and here the controls are unconstrained. The second one is that the flux of controls is zero (i.e., $\int_{\partial\Omega}\mathbf{y}\cdot\mathbf{n}=0$) \cite{gong_mateoes_stokes}. The zero flux condition comes naturally on control since we have an impressibility condition and we have only Dirichlet boundary condition in the PDE. So, to hold this zero flux condition, the authors choose only tangential control as the first choice and as the other choice the authors take zero flux condition itself as a constraint in the space. Also, it is observed in many Navier-Stokes Dirichlet control problem that the authors use either tangential control or the zero flux condition on the control, for e.g., one can see \cite{Fursikov,Gunzburger}. This zero flux condition on the control reduces the regularity of the control discussed in \cite{gong_mateoes_stokes}. To overcome this difficulty we introduce the Stokes equation with outflow condition and the control acts on the Dirichlet boundary only. Hence our control is more general and it has both the tangential and normal components. Also, we have introduced constraints in the control. Due to these constraints in the control, the optimal control satisfies a simplified Signorini problem:
 	\begin{subequations}\label{signorini_L2_f}
 	\begin{align}
 		-\rho\Delta \mathbf{y}&=\mathbf{0}\quad \text{in}\quad \Omega,\\
 		\mathbf{y}&=\mathbf{0}\quad \text{on}\quad \Gamma_D\cup \Gamma_N,\\
 		\mathbf{y}_a \leq \bm{\gamma}_0(\mathbf{y})&\leq \mathbf{y}_b\; \text{ a.e. on } \Gamma_C,
 	\end{align}
 	further the following holds for almost every $x\in \Gamma_C$:
 	\begin{align}
 		\text{if}\;  \mathbf{y}_a<\mathbf{y}(x) <\mathbf{y}_b \quad \text{then} \quad \big(\bm{\mu}(\mathbf{y})\big)(x)&=\mathbf{0},\label{sig_conditions1}\\
 		\text{if}\; \mathbf{y}_a\leq \mathbf{y}(x)<\mathbf{y}_b \quad \text{then} \quad \big(\bm{\mu}(\mathbf{y})\big)(x)&\geq \mathbf{0},\\
 		\text{if}\;  \mathbf{y}_a<\mathbf{y}(x)\leq \mathbf{y}_b \quad \text{then} \quad \big(\bm{\mu}(\mathbf{y})\big)(x)&\leq \mathbf{0}, \label{signorini_L2_e}
 	\end{align}
 \end{subequations}
\noindent
where the contact stress $\bm{\mu}(\mathbf{y})=\rho\frac{\partial \mathbf{y}}{\partial \mathbf{n}}-\frac{\partial \bm{\phi}}{\partial
	\mathbf{n}}-r\mathbf{n},$  $(\bm{\phi},r)$ is the adjoint variable and $\mathbf{y}_a$ , $\mathbf{y}_b$ are vectors in $\mathbb{R}^2$ and $\Gamma_C$, $\Gamma_D$ and $\Gamma_N$ are subsets of $\partial \Omega.$ As a result of this inequality in the contact boundary $\Gamma_C$, if we apply the standard error analysis for  $\|\nabla(\mathbf{y}-\mathbf{y}_h)\|_{0,\Omega},$ we only achieve $\frac{1}{2}+\frac{\delta}{2}$ $(\delta>0)$ order of convergence, when $\mathbf{y}\in \mathbf{H}^{\frac{3}{2}+\delta}(\Omega).$ However, this is not the optimal rate of convergence. Using the ideas in \cite{Hild:2015}, we have derived in Theorem \ref{energy_estm_control} that the control error $\|\nabla(\mathbf{y}-\mathbf{y}_h)\|_{0,\Omega}$ has $\frac{1}{2}+\delta$ order of convergence, which is optimal. Even if the control in the less regular space i.e., $\mathbf{y}\in \mathbf{H}^{\tau}(\Omega)$ with $1<\tau\leq 3/2$ we derive an optimal order of convergence up to the regularity. To this end, we have used a two-level finite element method for the Stokes problem where the velocity and the control are approximated by a piecewise linear finite element space on a fine mesh and the pressure is approximated by a piecewise constant space on a coarse mesh. This way the point-wise control constraints are well respected. The theoretical results are corroborated by a variety of numerical tests

\par

The rest of the article is structured as follows. In Section \ref{sec:ModelProblem}, we prove the existence and uniqueness of the solution of the optimal control problem and derive the continuous optimality system. In Section \ref{Discrete Problems}, we discuss about the discrete optimal control problem. We derive \textit{a priori} error estimates in Section \ref{Sto:DC_error_analysis}. Section
\ref{sec:Numerics} is devoted to the numerical experiments.

\section{Continuous Control Problem}\label{sec:ModelProblem}
In this section we briefly discuss the precise formulation of the optimization problem under consideration. Furthermore, we recall theoretical results on existence, uniqueness, and regularity of optimal solutions as well as optimality conditions. Before going to the analysis we need the following definitions:
\subsection{Notations}\label{notation}
Let $\Omega$ be a bounded convex polygonal domain in $\mathbb{R}^2$. Any function and space in bold notation can be understood in the vector form e.g.,
$\mathbf{x}:=(x_{1},x_{2}),$  $\mathbf{L}^2(\Omega):=[L^2(\Omega)]^2$and $\mathbf{H}^1(\Omega):=[H^1(\Omega)]^2.$ The norm and inner product on those spaces are defined component wise. The norm in the $\mathbf{L}^2(\Omega)$ space is denoted by $\|\cdot\|_{0,\Omega}.$ Also, the norm on the Sobolev space $\mathbf{H}^k(\Omega)$ is denoted by $\|\cdot\|_{k,\Omega} (k> 0),$ see \cite{Ciarlet:1978:FEM}. The trace of a vector valued function $\mathbf{x}\in \mathbf{H}^1(\Omega)$ is defined to be $\bm{\gamma}_{0}(\mathbf{x}):=(\gamma_{0}(x_1),\gamma_{0}(x_2)),$ where $\gamma_0:H^1(\Omega)\rightarrow L^2(\Gamma)$ is the
trace operator whose range is $H^{1/2}(\Gamma)$. Let  $\mathbf{x}$ and $\mathbf{y}$ are two functions, we say that  $\mathbf{x}\leq \mathbf{y}$ iff $x_1\leq y_1$ and $x_2\leq y_2$ almost everywhere in $\Omega.$

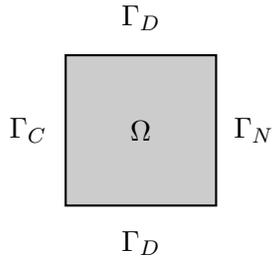
\begin{figure}[h!]
		\begin{center}
			\begin{tikzpicture}
			\fill[color=black!40, semitransparent]
			(0,0) -- (2,0) --(2,2) --(0,2)-- cycle;
			\draw[black,thick] (0,0) -- (2,0) --(2,2) --(0,2)  -- cycle;
			\node[black,thick] () at (-.5,1){$\Gamma_C$};
			\node[black,thick] () at (2.5,1){$\Gamma_N$};
			\node[black,thick] () at (1,2.5){$\Gamma_D$};
			\node[black,thick] () at (1,-.5){$\Gamma_D$};
			\node[black,thick] () at (1,1){$\Omega$};
			\end{tikzpicture}
			\caption{ The domain $\Omega.$}
			\label{example of domain}
		\end{center}
	\end{figure}

Before going to the optimal control problem first, we will discuss about the Stokes problem defined in \eqref{stokes_eq_intro}. Here we will describe the problem more precisely. We have the following Stokes problem with mixed boundary conditions:

\begin{subequations}\label{stokes_eq}
	\begin{align}
-\Delta \mathbf{u}+\nabla{p}&=\mathbf{f} \quad \text{in}\;\Omega,\label{EQI} \\
\nabla\cdot{\mathbf{u}}&=0 \quad \text{in}\; \Omega,\label{EQDI}\\
\bf{u}&=\mathbf{y} \quad \text{on}\; \Gamma_C,\label{EQCB}\\
\bf{u}&=\mathbf{0} \quad \text{on}\; \Gamma_D,\label{EQDB}\\
\frac{\partial \mathbf{u}}{\partial \mathbf{n}}-p \mathbf{n}&=\mathbf{0}\;\; \text{on} \;\; \Gamma_N. \label{EQNB}
\end{align}
\end{subequations}
Here $\Gamma_D, \;  \Gamma_C$ and $\Gamma_N$ are three non-overlapping open subsets of the boundary $\partial\Omega$ with $\partial\Omega=\Gamma_C\cup\bar \Gamma_D\cup\Gamma_N$, Figure \ref{example of domain} depicts an example of such a domain.
We assume that, $\Gamma_C$ be a straight line segment and one dimensional measure $|\Gamma_C|> 0$. The interior force $\mathbf{f}\in \mathbf{L}^2(\Omega)$.

\begin{remark}
	The choice of the domain $\Omega$ in this problem is very specific because of the mixed boundary conditions in the problem \ref{stokes_eq}. It is known that the Neumann-Dirichlet transition points always impair the regularity of the solution. The singularity of those transition points depends on the data as well as the interior angle at that point. Here, we have taken the angle between $\Gamma_D$ and $\Gamma_N$ is always $\pi/2$ and $\Gamma_C$ is a straight line segment so that, we can get a regular solution.
\end{remark}

\par
\noindent
Define the test and trial space $\mathbf{V}$ by $$\mathbf{V} :=\mathbf{H}_{D\cup C}^1(\Omega)=\{\mathbf{x} \in \mathbf{H}^1(\Omega): \bm{\gamma}_{0}(\mathbf{x})=\mathbf{0}\text{ on } \Gamma_D\cup \Gamma_C\}.$$
We choose controls from the following space:
$$\mathbf{Q}:=\{\mathbf{x} \in \mathbf{H}^1(\Omega): \bm{\gamma}_{0}(\mathbf{x})=\mathbf{0}\text{ on } \Gamma_D\cup \Gamma_N\}.$$
For given $\mathbf{y}\in \mathbf{Q}$, the weak formulation of \eqref{stokes_eq} is as follows: find $(\mathbf{w},p)\in \mathbf{V} \times L^2(\Omega)$ such that
\begin{subequations}\label{State_EQ}
\begin{align}
\mathbf{u}&=\mathbf{w}+\mathbf{y},\\
a(\mathbf{w},\mathbf{z})+b(\mathbf{z},p) &=( \mathbf{f},\mathbf{z})-a(\mathbf{y},\mathbf{z}) \;\;\;
{\rm for~all}\;\mathbf{z} \in \mathbf{V},\label{velocty_eqn} \\
b(\mathbf{w},q)&=-b(\mathbf{y},q)\;\quad {\rm for~all}\; q \in L^2(\Omega),\label{div_eqn}
\end{align}
\end{subequations}
where
 $a(\mathbf{w},\mathbf{z})=  \int_{\Omega} \nabla{\mathbf{w}}:\nabla{\mathbf{z}} \dx$ , $b(\mathbf{z},p)= - \int_{\Omega} p\nabla\cdot{\mathbf{z}} \dx,$ and the matrix product $A:B := \sum_{i,j=1}^{n}a_{ij}b_{ij}$ when $A=(a_{ij})_{1\leq i,j \leq n}$ and $B=(b_{ij})_{1\leq i,j \leq n}$
with $(\cdot,\cdot)$ denotes the $\mathbf{L}^2(\Omega).$

The Babuska-Brezzi theorem \cite{BScott:2008:FEM,Ciarlet:1978:FEM} ensures the existence and uniqueness of the solution of \eqref{State_EQ}.
We define the solution map $S$ by $S(\mathbf{f},\mathbf{x}):=\mathbf{u}$, where $\mathbf{x}\in \mathbf{Q}$ be given and $\mathbf{u}=\mathbf{w}+\mathbf{x}$ solves \eqref{State_EQ}.
\par
\noindent Here we consider the energy cost functional $J:\mathbf{H}^1(\Omega)\times \mathbf{Q} \rightarrow \R$, which is defined in the equation \eqref{min:j}. There the constant $\rho>0$ is the regularizing parameter and $\mathbf{u_{d}} \in \mathbf{L}^2(\Omega)$ is a given target function. We seek control from the following constrained set:
$$\mathbf{Q}_{ad}:=\{\mathbf{x}\in \mathbf{Q}: \mathbf{y}_a\leq\bm{\gamma}_{0}(\mathbf{x})\leq \mathbf{y}_b\text{ a.e. on } \Gamma_C\},$$
where $\mathbf{y}_a=(y^1_a,y^2_a)$, $\mathbf{y}_b=(y^1_b,y^2_b)\in\R^2$ satisfying $y^1_a <y^2_a$ and $y^1_b <y^2_b$. Furthermore, whenever the set $\Gamma_D$ is non empty, for compatibility we assume that $y^1_a,y^1_b\leq 0$ and $y^2_a,y^2_b\geq0$ in order that, the control set $\mathbf{Q}_{ad}$ is nonempty.

\smallskip
\subsection{The model problem} Find $(\mathbf{u},\mathbf{y})\in \mathbf{H}_{D}^1(\Omega) \times \mathbf{Q}_{ad}$ such that
\begin{align}\label{eq:min-J1}
J(\mathbf{u},\mathbf{y})=\min_{(\mathbf{v},\mathbf{x})\in \mathbf{H}_{D}^1(\Omega)\times \mathbf{Q}_{ad}} J(\mathbf{v},\mathbf{x}),
\end{align}
subject to the condition that $\mathbf{v}=S(\mathbf{f},\mathbf{x})$.

\medskip
\par
\noindent The reduced cost functional $j:\mathbf{Q}\rightarrow \R$ defined as
\begin{align}\label{eq:j}
j(\mathbf{x}):=\frac{1}{2}\|S(\mathbf{f},\mathbf{x})-\mathbf{u_d}\|_{0,\Omega}^2+\frac{\rho}{2}\|\nabla
\mathbf{x}\|_{0,\Omega}^2,\quad \mathbf{x}\in \mathbf{Q}_{ad},\quad\rho>0.
\end{align}
Differentiating the reduced cost functional $j$, we obtain
\begin{align*}
j\,'(\mathbf{y})(\mathbf{x})=\lim_{t\rightarrow
	0}\frac{j(\mathbf{y}+t\mathbf{x})-j(\mathbf{y})}{t}=(S(\mathbf{f},\mathbf{y})-\mathbf{u_d}, S(\mathbf{0},\mathbf{x}))+\rho\,a(\mathbf{y},\mathbf{x}).
\end{align*}

\begin{theorem} [\bf Existance and uniqueness of the solution]
	There exists a unique solution of the control problem \eqref{eq:min-J1}.
\end{theorem}
\begin{proof}
 It is clear that the cost functional $J$ is non negative. Set,  $$m=\inf_{(\mathbf{v},\mathbf{x})\in \mathbf{H}_{D}^1(\Omega)\times \mathbf{Q}_{ad}} J(\mathbf{v},\mathbf{x}).$$ Then there exists a minimizing sequence $(\mathbf{u}_n,\mathbf{y}_n)$ such that $J(\mathbf{u}_n,\mathbf{y}_n)$ converges to $m$ and $\mathbf{u}_n=S(\mathbf{f},\mathbf{y}_n).$ Since, the sequence $J(\mathbf{u}_n,\mathbf{y}_n)$ convergent we can say that the sequences $\|\mathbf{u}_n-\mathbf{u_d}\|_{0,\Omega}$ and  $\|\nabla
 \mathbf{y}_n\|_{0,\Omega}$ are also convergent and hence bounded. Now, $\mathbf{y}_n\in \mathbf{Q}_{ad}$ by using the Poincar\'e inequality we can conclude that the sequence $\mathbf{y}_n$ is bounded in $\mathbf{Q}.$ Then there
 exists a subsequence of $\mathbf{y}_n$, still indexed by $n$ to simplify the notation, and a function $\mathbf{y}$,
 such that $\mathbf{y}_n$ converges to $\mathbf{y}$ weakly in $\mathbf{Q}.$ It is clear that the set $\mathbf{Q}_{ad}$ is closed and convex so it is weakly closed. Hence, $\mathbf{y}\in \mathbf{Q}_{ad}.$ \textit{A priori} estimate of the Stokes problem \eqref{State_EQ} from \cite{F_Boyer:Book}, we get
 \begin{align}\label{stability:Estm}
 	\|\nabla \mathbf{w}_n\|_{0,\Omega}+\|p_n\|_{0,\Omega}\leq C (\norm{\mathbf{f}}_{0,\Omega}+\norm{\nabla\mathbf{y}_n}_{0,\Omega}),
 \end{align}
  where $\mathbf{u}_n=\mathbf{w}_n+\mathbf{y}_n.$ Since, the sequence $\mathbf{y}_n$ is bounded in $\mathbf{Q}$ and using \eqref{stability:Estm} we can conclude that the sequence $\mathbf{w}_n$ is bounded in $\mathbf{H}^1_{0}(\Omega).$ So, we can extract a subsequence of it and name it by $\mathbf{w}_n$ and it converges to $\mathbf{w}.$ Thus we can extract a subsequence of $p_n$ (call it $p_n$) correspond to the subsequence of $\mathbf{w}_n$ such that, $p_n$ weakly converges to $p$ in $L^2(\Omega).$ Now we need to show that $\mathbf{w}$ is the corresponding candidate for the control $\mathbf{y}.$
 We have
 $$\int_{\Omega}\nabla \mathbf{w}_n : \nabla \mathbf{v}+\int_{\Omega} p_n \nabla \cdot \mathbf{v}=\int_{\Omega} \mathbf{f}\cdot\mathbf{v}-\int_{\Omega}\nabla \mathbf{y}_n : \nabla \mathbf{v}\quad {\rm for~all}\;\mathbf{v} \in \mathbf{V}.$$
 Using the above weak convergences, we can conclude that
 $$\int_{\Omega}\nabla \mathbf{w}: \nabla \mathbf{v}+\int_{\Omega} p \nabla \cdot \mathbf{v}=\int_{\Omega} \mathbf{f}\cdot\mathbf{v}-\int_{\Omega}\nabla \mathbf{y}: \nabla \mathbf{v}\quad {\rm for~all}\;\mathbf{v} \in \mathbf{V}.$$
 Hence, $\mathbf{u}=\mathbf{w}+\mathbf{y}$ be the corresponding state for the control $\mathbf{y}.$
  The weak lower semi continuity of the norm gives, $\norm{\nabla \mathbf{y}}_{0,\Omega}\leq \liminf_{n\rightarrow \infty} \norm{\nabla \mathbf{y}_n}_{0,\Omega}.$ Using the above weak convergences of $\mathbf{w}_n$ and $\mathbf{y}_n$ we can conclude that both the sequences converges strongly in $\mathbf{L}^2(\Omega).$ Thus, $\mathbf{u}_n$ converges to $\mathbf{u}$ strongly in $\mathbf{L}^2(\Omega).$ Hence, we have
 $$J(\mathbf{u},\mathbf{y})\leq  \lim_{n\rightarrow \infty}\frac{1}{2} \|\mathbf{u}_n-\mathbf{u_d}\|_{0,\Omega}^2+\frac{\rho}{2}\liminf_{n\rightarrow \infty} \norm{\nabla \mathbf{y}_n}_{0,\Omega}^2=m.$$
  This proves the existence of a control $\mathbf{y}$ such that $J(\mathbf{u},\mathbf{y})=m.$  The uniqueness of the solution follows from the strict convexity of the cost functional.
\end{proof}

\noindent
The following proposition establishes the first order optimality system:

\begin{proposition}[\bf Continuous Optimality System]\label{prop:Wellposed-C}
	The state, adjoint state, and control $((\mathbf{u},p),(\bm{\phi},r),\mathbf{y})\in (\mathbf{H}_{D}^1(\Omega)\times L^2(\Omega))\times (\mathbf{V}\times L^2(\Omega)) \times \mathbf{Q}_{ad}$ satisfy the optimality system
	\begin{subequations}\label{conti_optimality_system}
	\begin{align}
	\mathbf{u}&=\mathbf{w}+\mathbf{y},\quad \mathbf{w}\in \mathbf{V},\label{eq:state_0a}\\
	a(\mathbf{w},\mathbf{z})+b(\mathbf{z},p) &=( \mathbf{f},\mathbf{z})-a(\mathbf{y},\mathbf{z}) \;\;\;
\forall\;\mathbf{z} \in \mathbf{V},\label{eq:state_a}\\ 
	b(\mathbf{u},q)&=0\;\quad \forall\; q \in  L^2(\Omega),\label{eq:state_b}\\
	a(\mathbf{z},\bm{\phi})-b(\mathbf{z},r) &=( \mathbf{u-u_d},\mathbf{z}) \;\;\;
	\forall\;\mathbf{z} \in \mathbf{V},\label{eq:adjstate_a}\\ 
	b(\bm{\phi},q)&=0\;\quad \forall\; q\in  L^2(\Omega),\label{eq:adjstate_b}\\
	\rho\, a(\mathbf{y},\mathbf{x}-\mathbf{y})&\geq a(\mathbf{x}-\mathbf{y},\bm{\phi})-b(\mathbf{x}-\mathbf{y},r)-(\mathbf{u-u_d},\mathbf{x}-\mathbf{y})\quad\forall \mathbf{x}\in
	\mathbf{Q}_{ad}.\label{eq:VI}
	\end{align}
\end{subequations}
\end{proposition}
\begin{proof} The equations from \eqref{eq:state_0a} to \eqref{eq:adjstate_b} are optimal state and adjoint state equations. Only we need to prove the last inequality \eqref{eq:VI}. The first order necessary optimality
	conditions yields
	\begin{align*}
	(\mathbf{u-u_d},S\mathbf{(0,x-y)})+\rho \,a\mathbf{(y,x-y)}\geq 0\quad \forall \mathbf{x}\in \mathbf{Q}_{ad}.
	\end{align*}
	The solution of the adjoint problem is defined by
	\begin{align*}
	a(\mathbf{z},\bm{\phi})-b(\mathbf{z},r) &=( \mathbf{u-u_d},\mathbf{z}) \;\;\;
	{\rm for~all}\;\mathbf{z} \in \mathbf{V},\\ \nonumber
	b(\bm{\phi},q)&=0\;\quad {\rm for~all}\; q \in L^2(\Omega)
	\end{align*}
	Since $S(\mathbf{0},\mathbf{x}-\mathbf{y})-(\mathbf{x}-\mathbf{y}) \in \mathbf{V}$ and $a(S(0,\mathbf{x}-\mathbf{y}),\bm{\phi})=0$, we obtain
	\begin{align*}
	(\mathbf{u-u_d},S\mathbf{(0,x-y)})&=(\mathbf{u-u_d},S\mathbf{(0,x-y)-(x-y)})+(\mathbf{u-u_d},\mathbf{x-y})\\
	&=(\mathbf{u-u_d},\mathbf{x}-\mathbf{y})-a(\mathbf{x}-\mathbf{y},\bm{\phi})+b(\mathbf{x-y,r}).
	\end{align*}
	This completes the proof.
\end{proof}

\begin{remark}[\bf Control Satisfies Signorini Problem]	
The optimal control $\mathbf{y}$ is the weak solutions of the  simplified Signorini problem \eqref{signorini_L2_f} defined in the introduction.
\end{remark}
\begin{remark}[\bf Regularity of Signorini Problem]\label{Rem_reg_signorini}
The numerical analysis of any finite element method applied to the Signorini problem \eqref{signorini_L2_f} requires
the knowledge of the regularity of the solution $\mathbf{y}$. Since the work by Moussaoui and Khodja (see \cite{Moussaoui_et_al}),
it is admitted that the Signorini condition may generate some singular behavior at the neighborhood of $\Gamma_C$. There are many factors that affect the regularity of the solution to the Signorini problem. Some of those factors are the regularity of the data, the mixed boundary conditions
(e.g., Neumann-Dirichlet transitions), the corners in polygonal domains and the Signorini condition, which generates singularities at contact-noncontact transition points. In our model problem we assume $\Gamma_{C}$ be a straight line segment. Let $\mathbf{p}$ be a contact-noncontact transition point in the interior of  $\Gamma_{C},$ then the solution of Signorini problem \eqref{signorini_L2_f} $\mathbf{y} \in \mathbf{H}^{\tau}(V_\mathbf{p})$ with $\tau < \frac{5}{2}$ and $V_\mathbf{p}$ be an open neighbourhood of $\mathbf{p}$ (see \cite[subsection 2.3]{Belgacem:2013}, \cite[section 2]{Belhachmi:2003} and \cite{Moussaoui_et_al}). Let $\mathbf{p}\in \bar{\Gamma}_{C}\cap \bar{\Gamma}_{D}$ and $V_\mathbf{p}$ be a neighbourhood of $\mathbf{p}$ in $\Omega$ such that $\mathbf{y}$ vanishes on $V_\mathbf{p}\cap\Gamma_{C}$ then the elliptic regularity theory on convex domain yields $\mathbf{y}\in \mathbf{H}^{2}(V_\mathbf{p})$(see \cite[subsection 2.3]{Belgacem:2013} and \cite{Orlt:1995}).
  Now if $\mathbf{y}$ does not vanish on $V_\mathbf{p}\cap\Gamma_{C},$ then $\mathbf{p}$ be a contact-noncontact type transition point and hence $\mathbf{y} \in \mathbf{H}^{\tau}(V_\mathbf{p})$ with $\tau < 5/2$ (see \cite[subsection 2.3]{Belgacem:2013} and \cite{Orlt:1995}). The best we can expect is to obtain $\mathbf{y} \in \mathbf{H}^{\tau}(V_{\Gamma_{C}})$ with $\tau \leq 2$ and $V_{\Gamma_{C}}$ is an open neighbourhood of $\Gamma_{C}$ (see \cite{Moussaoui_et_al,Belgacem:2013}).
\end{remark}

\begin{remark}[\bf Regularity of the adjoint state variables]\label{Rem_reg_adj_state} The strong form of the adjoint state is the following:
	\begin{subequations}\label{adj_state_eq}
		\begin{align*}
		-\Delta \bm{\phi}+\nabla{r}&=\mathbf{u}-\mathbf{u}_d \quad \text{in}\;\Omega, \\
		\nabla\cdot{\bm{\phi}}&=0 \quad \text{in}\; \Omega,\\
		\bm{\phi}&=\mathbf{0} \quad \text{on}\; \Gamma_D\cup \Gamma_C,\\
		\frac{\partial \bm{\phi}}{\partial \mathbf{n}}-r \mathbf{n}&=\mathbf{0}\;\; \text{on} \;\; \Gamma_N.
		\end{align*}
	\end{subequations}
Since our concerned domain has $\pi/2$ angle at Neumann-Dirichlet transition points and the given data are sufficiently regular so using Theorem A.1 of the paper \cite{Michal_et_al:2016} we can conclude that $\bm{\phi}\in \mathbf{H}^2(\Omega)$ and $r\in H^1(\Omega).$	
\end{remark}

\begin{remark}[\bf Regularity of the state variables]\label{Rem_reg_state}
We have seen form the Remark \eqref{Rem_reg_signorini} that the control variable can have regularity up to $\mathbf{H}^2(\Omega).$ So, if we assume $\mathbf{y}\in \mathbf{H}^2(\Omega)$ and since we have $\pi/2$ angle at each transition(Dirichlet-Dirichlet and Neumann-Dirichlet) points and load $\mathbf{f}\in \mathbf{L}^2(\Omega)$ then, from the equations \eqref{eq:state_a}-\eqref{eq:state_b} we can conclude that $\mathbf{w}\in \mathbf{H}^{\frac{3}{2}+\delta}(\Omega)$ and $p\in H^{\frac{1}{2}+\delta}(\Omega
)$ with $\delta>0$ (see, \cite{Orlt:1995}). Therefore, the velocity $\mathbf{u}\in \mathbf{H}^{\frac{3}{2}+\delta}(\Omega).$	
\end{remark}

\section{Discrete Control Problem}\label{Discrete Problems}
Let $\cT_H$ be a shape-regular triangulation of the domain $\Omega$ into triangles $K$ such that $\cup_{K\in \cT_H}K=\bar{\Omega}$ see \cite{BScott:2008:FEM,Ciarlet:1978:FEM}. Also let $\cT_h$ be a refinement of $\cT_H$ by connecting all the midpoints of $\cT_H.$ The collection of interior edges of $\cT_h$ is denoted by $\mathcal{E}_{h}^{i}.$ The collection of Dirichlet, Neumann and Contact boundary edges of $\cT_h$ are denoted by $\mathcal{E}_{h}^{b,D},$ $\mathcal{E}_{h}^{b,N}$ and $\mathcal{E}_{h}^{b,C}$ respectively. We define $\mathcal{E}_{h}=\mathcal{E}_{h}^{i}\cup\mathcal{E}_{h}^{b,D}\cup \mathcal{E}_{h}^{b,N}\cup\mathcal{E}_{h}^{b,C}.$ A typical triangle is denoted by $K$ and its
diameter by $h_T$. Set $h=\max_{K\in\cT_h} h_T$. The length of any edge $e\in \mathcal{E}_{h}$ will be denoted by $h_e$. The collection of all vertices of  $\cT_h$ is denoted by $\cV_h$. The set of vertices on $\overline{\Gamma}_D,$ $\Gamma_N$ and $\Gamma_C$ are denoted by $\cV_h^D,$ $\cV_h^N$ and $\cV_h^C.$ Define the discrete space for velocity $\mathbf{V}_h\subset \mathbf{V}$ by
$$\mathbf{V}_h:=\{\mathbf{v}_h\in \mathbf{V}: \mathbf{v}_h|_T\in \mathbf{P}_1(K)\;\; \forall K \in \mathcal{T}_h\},$$
and the discrete space for pressure is
$$M_H:=\{p_H\in L^2(\Omega): p_H|_T\in P_0(K)\;\; \forall K \in \mathcal{T}_H\},$$
and the discrete control space $\mathbf{Q}_h\subset \mathbf{Q}$ by
$$\mathbf{Q}_h:=\{\mathbf{x}_h\in \mathbf{Q}: \mathbf{x}_h|_T\in \mathbf{P}_1(K),\;\; \forall K \in \mathcal{T}_h\},$$
where $P_0(K)$ is the space of constant polynomials on $K$ and $\mathbf{P}_1(K)$ is the space of polynomials of degree $\leq 1$ on the triangle $K$. The approximation of the control set is as follows
$$\mathbf{Q}^h_{ad}:=\{\mathbf{x}_h\in \mathbf{Q}_h:\mathbf{y}_a\leq \mathbf{x}_h(z) \leq \mathbf{y}_b \text{ for all } z\in \cV_h^C\}.$$

It is easy to check that $\mathbf{Q}_{ad}^h\subset \mathbf{Q}_{ad}.$ Define the Lagrange interpolation $\mathcal{J}_h:\mathbf{C}(\bar{\Omega})\rightarrow \mathbf{Q}_h$ by $\mathcal{J}_h(\mathbf{v})=\sum_{z\in \cV_h} \mathbf{v}(z)\bm{\phi}_z,$ where $\bm{\phi}_z$ are basis functions of $\mathbf{Q}_h$ and $\mathbf{C}(\bar{\Omega})$ is the space of continuous functions on $\Omega$. From now on it will be assumed that $C>0$ be a constant independent of the mesh size $h$.
\begin{proposition}[\bf Discrete Optimality System]\label{prop:Wellposedness-discrete_optmality}
There exists a unique $((\mathbf{w}_h,p_H),(\bm{\phi}_h,r_H),\mathbf{y}_h)\in \big(\mathbf{V}_h\times M_H\big)\times \big(\mathbf{V}_h\times M_H\big) \times \mathbf{Q}^h_{ad}$ solves the following:
\begin{subequations}
\begin{align}
		\mathbf{u}_h&=\mathbf{w}_h+\mathbf{y}_h,\quad \mathbf{w}_h\in \mathbf{V}_h,\label{eq:d_state_0a}\\
		a(\mathbf{w}_h,\mathbf{z}_h)+b(\mathbf{z}_h,p_H) &=( \mathbf{f},\mathbf{z}_h)-a(\mathbf{y}_h,\mathbf{z}_h) \;\;\;
		{\rm for~all}\;\mathbf{z}_h \in \mathbf{V}_h,\label{eq:d_state_a}\\ 
		b(\mathbf{u}_h,q_H)&=0\;\quad {\rm for~all}\; q_H \in M_H ,\label{eq:d_state_b}\\
		a(\mathbf{z}_h,\bm{\phi}_h)-b(\mathbf{z}_h,r_H) &=( \mathbf{u}_h-\mathbf{u_d},\mathbf{z}_h) \;\;\;
		{\rm for~all}\;\mathbf{z}_h \in \mathbf{V}_h ,\label{eq:d_adjstate_a}\\ 
		b(\bm{\phi}_h,q_H)&=0\;\quad {\rm for~all}\; q_H \in M_H ,\label{eq:d_adjstate_b}\\
		\rho\, a(\mathbf{y}_h,\mathbf{x}_h-\mathbf{y}_h)&\geq a(\mathbf{x}_h-\mathbf{y}_h,\bm{\phi}_h)-b(\mathbf{x}_h-\mathbf{y}_h,r_H)\notag\\&-(\mathbf{u}_h-\mathbf{u_d},\mathbf{x}_h-\mathbf{y}_h)\;\quad {\rm for~all}\; \mathbf{x}_h\in
		\mathbf{Q}^h_{ad}.\label{eq:d_VI}
\end{align}
\end{subequations}
\end{proposition}
\begin{proof}
	Before going to prove the existence of the optimal solution we need to check the inf-sup stability of the bilinear form $b.$ The bilinear form $b$ is inf-sup stable for the pair $(\mathbf{V}_h,M_H)$ because of the existence of the Fortin operator (see \cite[Section 4.2.2]{Ern:2004:FEMbook}) $\bm{\pi}_h:\mathbf{V}\rightarrow \mathbf{V}_h$ defined by : For any vertex of the fine triangle which is a mid point of an edge $E$ of the coarse triangle, we define
	\begin{align*}
		\int_{E}(\bm{\pi}_h\mathbf{v}-\mathbf{v})=\mathbf{0}.
	\end{align*}
 Also, we define $\bm{\pi}_h\mathbf{v}$ at the vertices $\{a_1,a_2,a_3\}$ of the coarse triangle by using Scott-Zhang interpolation. Let $E_i$ be an edge containing the vertex $a_i.$ The set $\{ \phi_i^1,\phi_i^2\}$ denotes the restriction to $E_i$ of the local shape functions associated with the nodes lying in $E_i.$ Now consider the corresponding $L^2-$ dual basis of $\{ \phi_i^1,\phi_i^2\}$ is $\{ \psi_i^1,\psi_i^2\}$ such that
	\begin{align*}
		\int_{E_i}\psi^k_{i}\phi^l_{i}=\delta_{kl} \quad 1\leq k,l \leq 2.
	\end{align*}
Conventionally, set $\phi_i^1=\phi_i$ and $\psi_i^1=\psi_i$ for the node $a_i.$
	The nodal variables at vertices are defined by:
	$\bm{\pi}_h\mathbf{v}(a_i)=\int_{E_i} \psi_i\mathbf{v}.$
	Whenever $a_i$ is at the boundary and in the intersection of many edges,
	it is important to pick the one edge such that $E_i\subseteq \partial\Omega$.
	The map $\bm{\pi}_h$ defined above is a well defined operator. The $H^1$ stability of $\bm{\pi}_h$ can be derived using standard scaling argument see, \cite[Section 1.6.2]{Ern:2004:FEMbook}. The standard theory of optimal control problem \cite{trolzstch:2005:Book,Raymond:Notes} can be used to prove the existence and uniqueness of the solution.
\end{proof}
\section{Error Analysis}\label{Sto:DC_error_analysis}
In this section is devoted to \textit{a priori} error estimates. The convergence rate of the finite element approximation of the control problem depends on the regularity of the solution. It is clear from Remark \ref{Rem_reg_signorini}, \ref{Rem_reg_adj_state} and \ref{Rem_reg_state} that one can assume the solution $\mathbf{u}\in \mathbf{H}^{\frac{3}{2}+\delta}(\Omega),$ $p \in H^{\frac{1}{2}+\delta}(\Omega),$ $\bm{\phi}\in \mathbf{H}^{\frac{3}{2}+\delta}(\Omega),$  $r \in H^{\frac{1}{2}+\delta}(\Omega)$ and $\y \in \mathbf{H}^{\frac{3}{2}+\delta}(\Omega),$ where $ 0<\delta\leq 1/2.$ To derive \textit{a priori} error analysis, we introduce some projections as follows: Let $\mathbf{P}_{h}\mathbf{w} \in \mathbf{V}_{h},$ $\bar{\mathbf{P}}_{h}\bm{\phi} \in \mathbf{V}_{h},$ $R_{H}p \in M_{H}$ and $\bar{R}_{H}r \in M_{H}$ solve
\begin{subequations}
	\begin{align}
		a(\mathbf{P}_{h}\mathbf{w},\mathbf{z}_{h})+b(\mathbf{z}_{h},R_{H}p) &={\langle \mathbf{f},\mathbf{z}_{h} \rangle}- a(\mathbf{y},\mathbf{z}_{h}) \quad {\rm for~all}\;\mathbf{z}_{h} \in \mathbf{V}_{h},\label{eq:aux_state_a} \\
		b(\mathbf{P}_{h}\mathbf{w},q_{H})&=-b(\mathbf{y},q_{H})\;\quad {\rm for~all}\;q_{H} \in M_{H},\label{eq:aux_state_b}\\
		a_{h}(\mathbf{z}_{h},\bar{\mathbf{P}}_{h}\bm{\phi})-b_{h}(\mathbf{z}_{h},\bar{R}_{h}r) &={\langle \mathbf{u}-\mathbf{u}_d,\mathbf{z}_{h} \rangle}_{W} \;\quad{\rm for~all}\;\mathbf{z}_{h} \in \mathbf{V}_{h},\label{eq:aux_adjstate_a}\\
		b_{h}(\bar{\mathbf{P}}_{h}\bm{\phi},q_{H})&=0\;\quad{\rm for~all}\;q_{H} \in M_{H}\label{eq:aux_adjstate_b}.
	\end{align}
\end{subequations}
The following theorem is the first step to get error estimates.
\begin{theorem}[\bf Energy error estimate of control and $L^2$-estimate of velocity]\label{lem:Abstarct-error}
Let $\mathbf{y},\mathbf{u}$ be the continuous optimal control and state satisfy \eqref{conti_optimality_system}, $\mathbf{y}_h, \mathbf{u}_h$ be the discrete optimal control and state satisfy \eqref{prop:Wellposedness-discrete_optmality}. Then there holds
	\begin{align}\label{eq:Abstarct-error}
		\rho\|\nabla(\mathbf{y}-\mathbf{y}_h)\|_{0,\Omega}^2+&\|\mathbf{u}-\mathbf{u}_h\|_{0,\Omega}^2\leq
		C|a(\mathbf{y},\mathbf{y}-\mathbf{x}_h)-a(\mathbf{x}_h-\mathbf{y},\bm{\phi})+b(\mathbf{y}-\mathbf{x}_h,r)\nonumber\\
		&-(\mathbf{u}-\mathbf{u}_d,\mathbf{x}_h-\mathbf{y})|+C \left(\|\nabla(\mathbf{y}-\mathbf{x}_h)\|_{0,\Omega}^2+\|\nabla(\bm{\phi}-\bar{\mathbf{P}}_h\bm{\phi})\|_{0,\Omega}^2\right)\nonumber\\&+C\left(\|\mathbf{P}_h\mathbf{w}-\mathbf{w}\|_{0,\Omega}^2+\norm{r-\bar{R}_Hr}_{0,\Omega}^2+\|\mathbf{y}-\mathbf{x}_h\|_{0,\Omega}^2\right),
	\end{align}
	for all $\mathbf{x}_h\in \mathbf{Q}_{ad}^h.$
\end{theorem}

\begin{proof}
	A selection $\mathbf{x}=\mathbf{y}_h\in \mathbf{Q}^h_{ad}\subseteq \mathbf{Q}_{ad}$ in \eqref{eq:VI}, yields
	\begin{align}
		\rho\, a(\mathbf{y},\mathbf{y}_h-\mathbf{y})&\geq a(\mathbf{y}_h-\mathbf{y},\bm{\phi})-b(\mathbf{y}_h-\mathbf{y},r)-(\mathbf{u}-\mathbf{u}_d,\mathbf{y}_h-\mathbf{y}).\label{eq:EEP1}
	\end{align}
	Using \eqref{eq:d_VI}, we have
	\begin{align}
		\rho \,a(\mathbf{y}_h,\mathbf{y}-\mathbf{y}_h)&\geq -\rho \,a(\mathbf{y}_h,\mathbf{x}_h-\mathbf{y})+
		a(\mathbf{x}_h-\mathbf{y}_h,\bm{\phi}_h)-b(\mathbf{x}_h-\mathbf{y}_h,r_H)\notag\\&\quad -(\mathbf{u}_h-\mathbf{u}_d,\mathbf{x}_h-\mathbf{y}_h )\quad\forall
		\mathbf{x}_h\in \mathbf{Q}^h_{ad}. \label{eq:EEP2}
	\end{align}
	Adding the equations \eqref{eq:EEP1} and \eqref{eq:EEP2}, we find that for
	any $\mathbf{x}_h\in \mathbf{Q}^h_{ad}$
	\begin{align}
		\rho \,a(\mathbf{y}_h-\mathbf{y},\mathbf{y}-\mathbf{y}_h)&\geq  -\rho \,a(\mathbf{y}_h,\mathbf{x}_h-\mathbf{y})+
		a(\mathbf{x}_h-\mathbf{y}_h,\bm{\phi}_h)-b(\mathbf{x}_h-\mathbf{y}_h,r_H)\notag\\&-(\mathbf{u}_h-\mathbf{u}_d,\mathbf{x}_h-\mathbf{y}_h)+a(\mathbf{y}_h-\mathbf{y},\bm{\phi})-b(\mathbf{y}_h-\mathbf{y},r)-(\mathbf{u}-\mathbf{u}_d,\mathbf{y}_h-\mathbf{y})\notag\\
		&\geq  \rho \,a(\mathbf{y}-\mathbf{y}_h,\mathbf{x}_h-\mathbf{y})-\rho \,a(\mathbf{y},\mathbf{x}_h-\mathbf{y})+
		a(\mathbf{x}_h-\mathbf{y},\bm{\phi}_h-\bm{\phi})\notag\\&+a(\mathbf{x}_h-\mathbf{y},\bm{\phi})+a(\mathbf{y}-\mathbf{y}_h,\bm{\phi}_h-\bm{\phi})-b(\mathbf{x}_h-\mathbf{y}_h,r_H-r)-b(\mathbf{x}_h-\mathbf{y}_h,r)\notag\\&-b(\mathbf{y}-\mathbf{y}_h,r_H-r)-(\mathbf{u}_h-\mathbf{u}_d,\mathbf{x}_h-\mathbf{y})-(\mathbf{u}_h-\mathbf{u},\mathbf{y}-\mathbf{y}_h)\notag\\
		&\geq  \left(\rho \,a(\mathbf{y},\mathbf{y}-\mathbf{x}_h)+a(\mathbf{x}_h-\mathbf{y},\bm{\phi})+b(\mathbf{y}-\mathbf{x}_h,r)-(\mathbf{u}-\mathbf{u}_d,\mathbf{x}_h-\mathbf{y})\right)\notag\\&+\rho
		a(\mathbf{y}-\mathbf{y}_h, \mathbf{x}_h-\mathbf{y})+a(\mathbf{x}_h-\mathbf{y},\bm{\phi}_h-\bm{\phi})+a(\mathbf{y}-\mathbf{y}_h, \bm{\phi}_h-\bar{\mathbf{P}}_h\bm{\phi})\notag\\&+a(\mathbf{y}-\mathbf{y}_h, \bar{\mathbf{P}}_h\bm{\phi}-\bm{\phi})-b(\mathbf{x}_h-\mathbf{y},r_H-r)-b(\mathbf{y}-\mathbf{y}_h,r_H-r)\notag\\&-(\mathbf{u}_h-\mathbf{u},\mathbf{x}_h-\mathbf{y})+\norm{\mathbf{u}-\mathbf{u}_h}_{0,\Omega}^2+(\mathbf{u}_h-\mathbf{u},\mathbf{w}_h-\mathbf{w})\notag\\
		&\geq  \left(\rho \,a(\mathbf{y},\mathbf{y}-\mathbf{x}_h)+a(\mathbf{x}_h-\mathbf{y},\bm{\phi})+b(\mathbf{y}-\mathbf{x}_h,r)-(\mathbf{u}-\mathbf{u}_d,\mathbf{x}_h-\mathbf{y})\right)\notag\\&+\rho
		a(\mathbf{y}-\mathbf{y}_h, \mathbf{x}_h-\mathbf{y})+a(\mathbf{x}_h-\mathbf{y},\bm{\phi}_h-\bm{\phi})-b(\mathbf{y}_h-\mathbf{y},r-\bar{R}_Hr)\notag\\&+a(\mathbf{y}-\mathbf{y}_h, \bar{\mathbf{P}}_h\bm{\phi}-\bm{\phi})-b(\mathbf{x}_h-\mathbf{y},r_H-r)\notag\\&-(\mathbf{u}_h-\mathbf{u},\mathbf{x}_h-\mathbf{y})+\norm{\mathbf{u}-\mathbf{u}_h}_{0,\Omega}^2+(\mathbf{u}_h-\mathbf{u},\mathbf{P}_h\mathbf{w}-\mathbf{w}).\label{last_ineq}
	\end{align}
	Using Cauchy-Schwarz inequality and Young's inequality in the equation \eqref{last_ineq}, we obtain
	\begin{align}\label{Estm}
		\norm{\nabla(\mathbf{y}-\mathbf{y}_h)}_{0,\Omega}^2+\norm{\mathbf{u}-\mathbf{u}_h}_{0,\Omega}^2 \leq& C \left(\rho \,a(\mathbf{y},\mathbf{x}_h-\mathbf{y})+a(\mathbf{y}-\mathbf{x}_h,\bm{\phi})+b(\mathbf{x}_h-\mathbf{y},r)\right)\notag\\&-(\mathbf{u}-\mathbf{u}_d,\mathbf{y}-\mathbf{x}_h)+\norm{\nabla(\mathbf{y}-\mathbf{x}_h)}_{0,\Omega}^2+\norm{\mathbf{y}-\mathbf{x}_h}_{0,\Omega}^2\notag\\&+\norm{\nabla{(\bm{\phi}-\bm{\phi}_h)}}_{0,\Omega}^2+\norm{\nabla{(\bar{\mathbf{P}}_h\bm{\phi}-\bm{\phi})}}^2+\norm{r-\bar{R}_Hr}_{0,\Omega}^2\notag\\&+\norm{r-r_H}_{0,\Omega}^2+\norm{\mathbf{P}_h\mathbf{w}-\mathbf{w}}_{0,\Omega}^2.
	\end{align}
	We need to estimate the terms $\norm{\nabla{(\bm{\phi}-\bm{\phi}_h)}}_{0,\Omega}$ and $\norm{r-r_H}_{0,\Omega}.$  Introducing the projection we have
	\begin{align}\label{adj_proj}
		\norm{\nabla{(\bm{\phi}-\bm{\phi}_h)}}_{0,\Omega}\leq\norm{\nabla{(\bm{\phi}-\bar{\mathbf{P}}_h\bm{\phi})}}_{0,\Omega}+\norm{\nabla{(\bar{\mathbf{P}}_h\bm{\phi}-\bm{\phi}_h)}}_{0,\Omega}.
	\end{align}
	A subtraction of \eqref{eq:adjstate_a} from  \eqref{eq:aux_adjstate_a} yields
	\begin{align*}
		a(\mathbf{v}_h,\bar{\mathbf{P}}_h\bm{\phi}-\bm{\phi}_h)+b(\mathbf{v}_h,r_H-\bar{R}_Hr)=(\mathbf{u}-\mathbf{u}_h,\mathbf{v}_h).
	\end{align*}
	By taking $\mathbf{v}_h=\bar{\mathbf{P}}_h\bm{\phi}-\bm{\phi}_h$ in the above equation and using the fact that $b(\bar{\mathbf{P}}_h\bm{\phi}-\bm{\phi}_h,r_H-\bar{R}_Hr)=0,$ we get $
	a(\bar{\mathbf{P}}_h\bm{\phi}-\bm{\phi}_h,\bar{\mathbf{P}}_h\bm{\phi}-\bm{\phi}_h)=(\mathbf{u}-\mathbf{u}_h,\bar{\mathbf{P}}_h\bm{\phi}-\bm{\phi}_h).$
	Applying Cauchy-Schwarz inequality we find
	\begin{align}\label{aux_adj_est}
		\norm{\nabla{(\bar{\mathbf{P}}_h\bm{\phi}-\bm{\phi}_h)}}_{0,\Omega}\leq \norm{\mathbf{u}-\mathbf{u}_h}_{0,\Omega}.
	\end{align}
	Hence,
	\begin{align}\label{adj_est}
		\norm{\nabla{(\bm{\phi}-\bm{\phi}_h)}}_{0,\Omega}\leq\norm{\nabla{(\bm{\phi}-\bar{\mathbf{P}}_h\bm{\phi})}}_{0,\Omega}+\norm{\mathbf{u}-\mathbf{u}_h}_{0,\Omega}.
	\end{align}
	To estimate the term $\norm{r-r_H}_{0,\Omega},$ we introduce the projection
	\begin{align}\label{adj_pressure_proj}
		\norm{r-r_H}_{0,\Omega}\leq \norm{r-\bar{R}_Hr}_{0,\Omega}+\norm{\bar{R}_Hr-r_H}_{0,\Omega}.
	\end{align}
	Using the inf-sup condition, we have
	\begin{align}\label{thm345}
		\beta \lVert \bar{R}_Hr-r_H\rVert_{0,\Omega}
		&\leq  \lVert \mathbf{u}-\mathbf{u}_h\rVert_{0,\Omega}+\norm{\nabla{(\bar{\mathbf{P}}_h\bm{\phi}-\bm{\phi}_h)}}_{0,\Omega}
		\leq 2\lVert \mathbf{u}-\mathbf{u}_h\rVert_{0,\Omega}.
	\end{align}
	In the above we have used the equation \eqref{aux_adj_est}.
	Hence we have the following:
	\begin{align}\label{adj_pressure_estm}
		\norm{r-r_H}_{0,\Omega}\leq  \norm{r-\bar{R}_Hr}_{0,\Omega}+2\lVert \mathbf{u}-\mathbf{u}_h\rVert_{0,\Omega}.
	\end{align}
	Substituting \eqref{adj_est} and \eqref{adj_pressure_estm} in \eqref{Estm} we get the desired result.
\end{proof}
Before going to derive the estimates for the terms on the right hand side of the equation \eqref{eq:Abstarct-error} we need to introduce few notations. Let $K$ be a triangle which shares an edge with $\Gamma_C$,  define $$\mathcal{N}=\{x\in K \cap\Gamma_C:\; \mathbf{y}_a<\mathbf{y}(x)<\mathbf{y}_b\},$$ and
$$\mathcal{C}=\{x\in K \cap\Gamma_C:\; \mathbf{y}(x)=\mathbf{y}_a\}\cup\{x\in K \cap\Gamma_C:\; \mathbf{y}(x)=\mathbf{y}_b\}.$$
The sets $\mathcal{C}$ and $\mathcal{N}$ are measurable as the function $\mathbf{y}|_{\Gamma_C}$ is continuous on $\Gamma_C.$ Let $|\mathcal{C}|$ and $|\mathcal{N}|$ be their measures. Now we prove the following lemma, which will be useful in the error estimation of the control.
\begin{lemma} \label{Hilld_trick}
	There holds
	\begin{align*}
	|\rho \,a(\mathbf{y},\mathbf{x}_h-\mathbf{y})+a(\mathbf{y}-\mathbf{x}_h,\bm{\phi})+b(\mathbf{x}_h-\mathbf{y},r)-(\mathbf{u}-\mathbf{u}_d,&\mathbf{y}-\mathbf{x}_h)|\leq Ch^{1+2\delta}\big(\norm{\y}^2_{\frac{3}{2}+\delta,\Omega}\\&+\norm{\bm{\phi}}^2_{\frac{3}{2}+\delta,\Omega}+\norm{r}^2_{\frac{1}{2}+\delta,\Omega}\big).	
	\end{align*}
\end{lemma}
\begin{proof}
	A use of adjoint PDE (in Remark \ref{Rem_reg_adj_state}), equation \eqref{signorini_L2_f}, and integration by parts yields
	\begin{align}\label{signori_term}
	\rho\,a(\mathbf{y},\mathbf{x}_h-\mathbf{y})+a(\mathbf{y}-\mathbf{x}_h,\bm{\phi})+b(\mathbf{x}_h-\mathbf{y},r)-(\mathbf{u}&-\mathbf{u}_d,\mathbf{y}-\mathbf{x}_h)\nonumber\\&
	=\int_{\Gamma_C}\bm{\mu}(\mathbf{y})\;(\mathbf{x}_h-\mathbf{y})ds.
	\end{align}
Choose $\mathbf{x}_h=\mathcal{J}_h\mathbf{y} \in \mathbf{Q}_h,$ then the right hand side of \eqref{signori_term} reads as and equals
	$$\int_{\Gamma_C} \bm{\mu}(\mathbf{y}) (\mathcal{J}_h\mathbf{y}-\mathbf{y})ds=\sum_{K\in \mathcal{T}_{h}} \int_{K \cap\Gamma_C} \bm{\mu}(\mathbf{y}) (\mathcal{J}_h\mathbf{y}-\mathbf{y})ds .$$
	
	\noindent	
	Therefore it remains to estimate the following:
	\begin{equation}\label{eqn:estimate}
		\int_{K \cap\Gamma_C} \bm{\mu}(\mathbf{y}) (\mathcal{J}_h\mathbf{y}-\mathbf{y})ds \quad \text{for all}\; K\in\cT_h.
	\end{equation}
	\noindent
	Fix a triangle $K,$ sharing an edge with $\Gamma_{C}$. Denote the length of the edge $e=T\cap\Gamma_{C}$ by $h_e$. Clearly, $|\mathcal{C}|+|\mathcal{N}|=h_e.$ Now, if either $|\mathcal{C}|$ or $|\mathcal{N}|$ equals zero, then it is easy to see that the integral term in \eqref{eqn:estimate} vanishes. So we suppose that both $\mathcal{C}$ and $\mathcal{N}$ have positive measure in the following estimation of \eqref{eqn:estimate}. Now we will derive some estimates for the term \eqref{eqn:estimate}:\\
	\noindent
	\textit{Estimate of (\ref{eqn:estimate}) depending on $\mathcal{N}$:} Applying Cauchy-Schwarz inequality, and estimation in \eqref{lemma_result_1} in Lemma \ref{lemma_req} (see the subsection \ref{Appendix}), and a standard interpolation estimate yields
	\begin{align}
		\int_{K \cap\Gamma_C} \bm{\mu}(\mathbf{y}) (\mathcal{J}_h\mathbf{y}-\mathbf{y})ds &\leq \norm{\bm{\mu}(\mathbf{y})}_{0,K\cap \Gamma_C} \norm{\mathcal{J}_h\mathbf{y}-\mathbf{y}}_{0,K\cap\Gamma_C}\nonumber\\&\leq C \frac{1}{|\mathcal{N}|^\frac{1}{2}} h_e^{\frac{1}{2}+\delta} |\bm{\mu}(\mathbf{y})|_{\delta,K \cap\Gamma_C} h^{1+\delta} |\mathbf{y}'|_{\delta,K \cap\Gamma_C}\nonumber\\& \leq  C \frac{1}{|\mathcal{N}|^\frac{1}{2}} h_e^{\frac{3}{2}+2\delta} \big(|\bm{\mu}(\mathbf{y})|^2_{\delta,K \cap\Gamma_C} + |\mathbf{y}'|^2_{\delta,K \cap\Gamma_C}\big).\label{1}
	\end{align}
	\textit{Estimate of (\ref{eqn:estimate}) depending on $\mathcal{C}$:} Using interpolation error estimation of $\mathcal{J}_h$, and estimations \eqref{lemma_result_1} and \eqref{lemma_result_4} in
	Lemma \ref{lemma_req} (see the subsection \ref{Appendix}), we obtain
	\begin{align}
		\int_{K \cap\Gamma_C} \bm{\mu}(\mathbf{y}) (\mathcal{J}_h\mathbf{y}-\mathbf{y})ds &\leq \norm{\bm{\mu}(\mathbf{y})}_{0,K\cap \Gamma_C} \norm{\mathcal{J}_h\mathbf{y}-\mathbf{y}}_{0,K\cap\Gamma_C}\nonumber\\&\leq C \norm{\bm{\mu}(\mathbf{y})}_{0,K\cap \Gamma_C} h_e^\frac{1}{2} \norm{\mathbf{y}'}_{L^1(K \cap\Gamma_C)}\nonumber\\& \leq  C \frac{1}{|\mathcal{C}|^\frac{1}{2}} h_e^{\frac{3}{2}+2\delta} \big(|\bm{\mu}(\mathbf{y})|^2_{\delta,K \cap\Gamma_C} + |\mathbf{y}'|^2_{\delta,K \cap\Gamma_C}\big).\label{2}
	\end{align}
	It is clear that either $|\mathcal{N}|$ or $|\mathcal{C}|$ is greater than or equal to $h_e /2$. Now by choosing the compatible estimation \eqref{1} or \eqref{2}, we obtain
	\begin{equation*}
		\int_{K \cap\Gamma_C} \bm{\mu}(\mathbf{y}) (\mathcal{J}_h\mathbf{y}-\mathbf{y})ds\leq  C  h_e^{1+2\delta} \big(|\bm{\mu}(\mathbf{y})|^2_{\delta,K \cap\Gamma_C} + |\mathbf{y}'|^2_{\delta,K \cap\Gamma_C}\big).
	\end{equation*}
By summing over all the triangles sharing an edge with $\Gamma_{C}$ and applying trace results we obtain
	\begin{equation*}
		\int_{\Gamma_C} \bm{\mu}(\mathbf{y}) (\mathcal{J}_h\mathbf{y}-\mathbf{y})ds\leq  C h^{1+2\delta} \big(|\bm{\mu}(\mathbf{y})|^2_{\delta,\Gamma_C} + |\mathbf{y}'|^2_{\delta,\Gamma_C}\big)\leq Ch^{1+2\delta} \left(\norm{\mathbf{y}}^2_{\frac{3}{2}+\delta,\Omega}+\norm{\bm{\phi}}^2_{\frac{3}{2}+\delta,\Omega}+\norm{r}^2_{\frac{1}{2}+\delta,\Omega}\right).
	\end{equation*}
	This completes the proof.
\end{proof}
\begin{theorem}[\bf Energy error estimate of control and $L^2$-estimate of velocity]\label{energy_estm_control}
Let $\mathbf{y},\mathbf{u}$ be the continuous optimal control and state \eqref{conti_optimality_system}, $\mathbf{y}_h, \mathbf{u}_h$ be the discrete optimal control and state \eqref{prop:Wellposedness-discrete_optmality}. Then there holds
	\begin{align*}
		\rho\norm{\nabla(\mathbf{y}-\mathbf{y}_h)}_{0,\Omega}+\|\mathbf{u}-\mathbf{u}_h\|_{0,\Omega}\leq& C\big(h^{\frac{1}{2}+\delta}\norm{\y}_{\frac{3}{2}+\delta,\Omega}+h^{\frac{1}{2}+\delta}\norm{\bm{\phi}}_{\frac{3}{2}+\delta,\Omega}+h^{\frac{1}{2}+\delta}\norm{r}_{\frac{1}{2}+\delta,\Omega}\\&+h^{\frac{3}{2}+\delta}\norm{\y}_{\frac{3}{2}+\delta,\Omega}+h^{\frac{3}{2}+\delta}\norm{\mathbf{u}}_{\frac{3}{2}+\delta,\Omega}\big).
	\end{align*}
\end{theorem}

\begin{proof}
	From Theorem \ref{lem:Abstarct-error} we have,
	\begin{align}\label{L2 estimate}
		\rho\|\nabla(\mathbf{y}-\mathbf{y}_h)\|_{0,\Omega}^2+&\|\mathbf{u}-\mathbf{u}_h\|_{0,\Omega}^2\leq
		C|a(\mathbf{y},\mathbf{y}-\mathbf{x}_h)-a(\mathbf{x}_h-\mathbf{y},\bm{\phi})+b(\mathbf{y}-\mathbf{x}_h,r)\nonumber\\&-(\mathbf{u}-\mathbf{u}_d,\mathbf{x}_h-\mathbf{y})|
		+C \left(\|\nabla(\mathbf{y}-\mathbf{x}_h)\|_{0,\Omega}^2+\|\nabla(\bm{\phi}-\bar{\mathbf{P}}_h\bm{\phi})\|_{0,\Omega}^2\right)\nonumber\\&+C\left(\|\mathbf{P}_h\mathbf{w}-\mathbf{w}\|_{0,\Omega}^2+\norm{r-\bar{R}_Hr}_{0,\Omega}^2+\|\mathbf{y}-\mathbf{x}_h\|_{0,\Omega}^2\right),
	\end{align}
	for all $\mathbf{x}_h\in \mathbf{Q}_{ad}^h.$ The first term in the right hand side of \eqref{L2 estimate} has been estimated in Lemma \ref{Hilld_trick}. The estimate of $\|\nabla(\bm{\phi}-\bar{\mathbf{P}}_h\bm{\phi})\|_{0,\Omega},$ $\|\mathbf{P}_h\mathbf{w}-\mathbf{w}\|_{0,\Omega},$ and $\norm{r-\bar{R}_Hr}_{0,\Omega}$ follows from \cite{Ern:2004:FEMbook}. A selection, $\mathbf{x}_h=\mathcal{J}_h\mathbf{y}$ gives the required estimate of the second and the last term. Using all the estimates together we achieve the desired estimate. 	
\end{proof}
\begin{theorem}[\bf Energy error estimate of velocity]\label{energy estm velocity}
	Let $\mathbf{u}$ be the continuous optimal velocity satisfies \eqref{eq:state_0a} and $\mathbf{u}_h$ be the discrete optimal velocity satisfies \eqref{eq:d_state_0a}. Then there holds
	\begin{align*}
		\norm{\nabla(\mathbf{u}-\mathbf{u}_h)}_{0,\Omega}\leq& C\big(h^{\frac{1}{2}+\delta}\norm{\y}_{\frac{3}{2}+\delta,\Omega}+h^{\frac{1}{2}+\delta}\norm{\bm{\phi}}_{\frac{3}{2}+\delta,\Omega}+h^{\frac{1}{2}+\delta}\norm{r}_{\frac{1}{2}+\delta,\Omega}\\&+h^{\frac{3}{2}+\delta}\norm{\y}_{\frac{3}{2}+\delta,\Omega}+h^{\frac{3}{2}+\delta}\norm{\mathbf{u}}_{\frac{3}{2}+\delta,\Omega}\big).
	\end{align*}
\end{theorem}
\begin{proof} Splitting the state $\mathbf{u=w+y}$ and the discrete state $\mathbf{u}_h=\mathbf{w}_h+\mathbf{y}_h,$ we have
	\begin{align*}
		\norm{\nabla(\mathbf{u}-\mathbf{u}_h)}_{0,\Omega}\leq \norm{\nabla(\mathbf{w}-\mathbf{w}_h)}_{0,\Omega}+\norm{\nabla(\mathbf{y}-\mathbf{y}_h)}_{0,\Omega}.
	\end{align*}
 Introducing the projection in the first term of the above equation we obtain $$\norm{\nabla(\mathbf{w}-\mathbf{w}_h)}_{0,\Omega}\leq \norm{\nabla(\mathbf{w}-\mathbf{P}_h\mathbf{w})}_{0,\Omega}+\norm{\nabla(\mathbf{P}_h\mathbf{w}-\mathbf{w}_h)}_{0,\Omega}.$$ Subtraction of \eqref{eq:d_state_a} from \eqref{eq:aux_state_a} yields
	\begin{align}
		a(\mathbf{P}_{h}\mathbf{w}-\mathbf{w}_h,\mathbf{z}_{h})+b(\mathbf{z}_{h},R_{H}p-p_H) = a(\mathbf{y}_h-\mathbf{y},\mathbf{z}_{h}) \quad {\rm for~all}\;\mathbf{z}_{h} \in \mathbf{V}_{h}.
	\end{align}
	By taking $\mathbf{z}_h=\mathbf{P}_h\mathbf{w}-\mathbf{w}_h$ in the above equation and using the fact that $b(\mathbf{P}_h\mathbf{w}-\mathbf{w}_h,R_H p-p_H)=0 $ we get
	$\norm{\nabla(\mathbf{P}_{h}\mathbf{w}-\mathbf{w}_h)}_{0,\Omega}^2=a(\mathbf{y}_h-\mathbf{y},\mathbf{P}_h\mathbf{w}-\mathbf{w}_h).$
	Applying Cauchy-Schwarz inequality we find
	\begin{align}\label{aux_velocity}
		\norm{\nabla(\mathbf{P}_{h}\mathbf{w}-\mathbf{w}_h)}_{0,\Omega}\leq \norm{\nabla(\mathbf{y}_h-\mathbf{y})}_{0,\Omega}.
	\end{align}
	Hence,
	\begin{align*}\label{Estm.1}
		\norm{\nabla(\mathbf{u}-\mathbf{u}_h)}_{0,\Omega}\leq \norm{\nabla(\mathbf{w}-\mathbf{P}_h\mathbf{w})}_{0,\Omega}+2\norm{\nabla(\mathbf{y}-\mathbf{y}_h)}_{0,\Omega}.
	\end{align*}
	Using the estimate of $\norm{\nabla(\mathbf{y}-\mathbf{y}_h)}_{0,\Omega}$ from  Theorem \ref{lem:Abstarct-error} and estimate of  $\norm{\nabla(\mathbf{w}-\mathbf{P}_h\mathbf{w})}_{0,\Omega}$ from \cite{Ern:2004:FEMbook} we have the required result.
\end{proof}	

\begin{theorem}[\bf Error estimate of pressure]\label{estimate of pressure}
Let $p$ be the continuous optimal pressure satisfies \eqref{eq:state_a} and $p_H$ be the discrete optimal pressure satisfies \eqref{eq:d_state_a}. Then there holds
	\begin{align*}
		\norm{p-p_H}_{0,\Omega}\leq& C\big(h^{\frac{1}{2}+\delta}\norm{\y}_{\frac{3}{2}+\delta,\Omega}+h^{\frac{1}{2}+\delta}\norm{\bm{\phi}}_{\frac{3}{2}+\delta,\Omega}+h^{\frac{1}{2}+\delta}\norm{r}_{\frac{1}{2}+\delta,\Omega}\\&+h^{\frac{3}{2}+\delta}\norm{\y}_{\frac{3}{2}+\delta,\Omega}+h^{\frac{3}{2}+\delta}\norm{\mathbf{u}}_{\frac{3}{2}+\delta,\Omega}\big).
	\end{align*}
	
\end{theorem}

\begin{proof}
	Introducing the projection we have, $\norm{p-p_H}_{0,\Omega}\leq \norm{p-R_Hp}_{0,\Omega}+\norm{R_Hp-p_H}_{0,\Omega}.$ The estimate of $\norm{R_Hp-p_H}_{0,\Omega}$ follows from the following:
	\begin{align}
		\beta \lVert R_Hp-p_H\rVert_{0,\Omega}&\leq \sup_{\mathbf{v}_h\in \mathbf{V}_h}\frac{b(\mathbf{v}_h, R_Hp-p_H)}{\lVert \mathbf{v}_h \rVert_{1}}\nonumber\\
		&\leq \sup_{\mathbf{v}_h\in \mathbf{V}_h} \frac{a( \mathbf{y}_h-\mathbf{y},\mathbf{v}_h) -a(\mathbf{P}_h\mathbf{w}-\mathbf{w}_h,\mathbf{v}_h )}{\lVert \mathbf{v}_h\rVert_{1}}\nonumber\\
		&\leq  \lVert \nabla (\mathbf{y}-\mathbf{y}_h)\rVert_{0,\Omega}+\norm{\nabla{(\mathbf{P}_h\mathbf{w}-\mathbf{w}_h)}}_{0,\Omega}.
	\end{align}
	Using \eqref{aux_velocity} in the above equation we have
	$\norm{p-p_H}_{0,\Omega}\leq  \norm{p-R_Hp}_{0,\Omega}+2\lVert \nabla (\mathbf{y}-\mathbf{y}_h)\rVert_{0,\Omega}.$
	Using the estimates of $\norm{p-R_Hp}_{0,\Omega}$ from \cite{Ern:2004:FEMbook} and $\lVert \nabla (\mathbf{y}-\mathbf{y}_h)\rVert_{0,\Omega}$ from Theorem \ref{energy_estm_control} we get the desired result.
\end{proof}
\begin{theorem}[\bf Error estimate of adjoint velocity]\label{estimate of adj ver}
Let $\bm{\phi}$ be the continuous optimal adjoint velocity satisfies \eqref{eq:adjstate_a} and $\bm{\phi}_h$ be the discrete optimal adjoint velocity satisfies \eqref{eq:d_adjstate_a}. Then there holds
	\begin{align*}
		\norm{\nabla(\bm{\phi}-\bm{\phi}_h)}_{0,\Omega}\leq& C\big(h^{\frac{1}{2}+\delta}\norm{\y}_{\frac{3}{2}+\delta,\Omega}+h^{\frac{1}{2}+\delta}\norm{\bm{\phi}}_{\frac{3}{2}+\delta,\Omega}+h^{\frac{1}{2}+\delta}\norm{r}_{\frac{1}{2}+\delta,\Omega}\\&+h^{\frac{3}{2}+\delta}\norm{\y}_{\frac{3}{2}+\delta,\Omega}+h^{\frac{3}{2}+\delta}\norm{\mathbf{u}}_{\frac{3}{2}+\delta,\Omega}\big).
	\end{align*}
	
\end{theorem}

\begin{proof}
	Using the equations \eqref{adj_proj}-\eqref{adj_est} we get
	\begin{align*}
		\norm{\nabla{(\bm{\phi}-\bm{\phi}_h)}}_{0,\Omega}\leq\norm{\nabla{(\bm{\phi}-\bar{\mathbf{P}}_h\bm{\phi})}}_{0,\Omega}+\norm{\mathbf{u}-\mathbf{u}_h}_{0,\Omega}.
	\end{align*}
	Applying Theorem \ref{energy_estm_control} and using the estimate of $\norm{\nabla{(\bm{\phi}-\bar{\mathbf{P}}_h\bm{\phi})}}_{0,\Omega}$ from \cite{Ern:2004:FEMbook}  we achieve the desired result.
\end{proof}

\begin{theorem}[\bf Error estimate of adjoint pressure]\label{estimate of adj pressure}
Let $r$ be the continuous optimal adjoint pressure satisfies \eqref{eq:adjstate_a} and $r_H$ be the discrete optimal adjoint pressure satisfies \eqref{eq:d_adjstate_a}. Then there holds
	\begin{align*}
		\norm{r-r_H}_{0,\Omega}\leq& C\big(h^{\frac{1}{2}+\delta}\norm{\y}_{\frac{3}{2}+\delta,\Omega}+h^{\frac{1}{2}+\delta}\norm{\bm{\phi}}_{\frac{3}{2}+\delta,\Omega}+h^{\frac{1}{2}+\delta}\norm{r}_{\frac{1}{2}+\delta,\Omega}\\&+h^{\frac{3}{2}+\delta}\norm{\y}_{\frac{3}{2}+\delta,\Omega}+h^{\frac{3}{2}+\delta}\norm{\mathbf{u}}_{\frac{3}{2}+\delta,\Omega}\big).
	\end{align*}
	
\end{theorem}

\begin{proof}
	Similar to the proof of Theorem \ref{estimate of pressure}. 	
\end{proof}

\begin{remark}\label{estm_less_reg}
	If the domain is not so smooth, e.g. the angle at transition (Dirichlet-Dirichlet and Neumann-Dirichlet) point is greater than $\pi/2$ or some bad polygonal structure,
	then there is a possibility that the solution could be less regular i.e., $\mathbf{u}\in \mathbf{H}^{\frac{3}{2}-\delta}(\Omega),$ $\bm{\phi}\in \mathbf{H}^{\frac{3}{2}-\delta}(\Omega),$ $r \in H^{\frac{1}{2}-\delta}(\Omega),$ and $\y \in \mathbf{H}^{\frac{3}{2}-\delta}(\Omega),$ where $ 0<\delta<1/2.$ Then all the above \textit{a priori} estimates hold true except Lemma \ref{Hilld_trick}. It is clear that if the solutions have the above regularity then \eqref{signori_term} is not true because the right hand side of \eqref{signori_term} does not make sense. So, to estimate the term
	\begin{align}\label{signorini_term}
		|\left(\rho \,a(\mathbf{y},\mathbf{x}_h-\mathbf{y})+a(\mathbf{y}-\mathbf{x}_h,\bm{\phi})+b(\mathbf{x}_h-\mathbf{y},r)-(\mathbf{u}-\mathbf{u}_d,\mathbf{y}-\mathbf{x}_h)\right)|
	\end{align}
	we use the following idea:
	\begin{align}
		\rho \,a(\mathbf{y},\mathbf{x}_h-\mathbf{y})+a(\mathbf{y}-\mathbf{x}_h,\bm{\phi})+b(\mathbf{x}_h-\mathbf{y},r)-(\mathbf{u}&-\mathbf{u}_d,\mathbf{y}-\mathbf{x}_h)\nonumber\\&
		=\langle\rho\frac{\partial \mathbf{y}}{\partial \mathbf{n}}-\frac{\partial \bm{\phi}}{\partial
			\mathbf{n}}-r\mathbf{n},\mathbf{x}_h-\mathbf{y}\rangle_{\delta,\Gamma_C}\nonumber\\
		&\leq \norm{\bm{\mu}(\mathbf{y})}_{H^{\delta}(\Gamma_C)'} \norm{\mathbf{x}_h-\mathbf{y}}_{\delta,\Gamma_C}.\label{signorini_less_regu}
	\end{align}
	Choosing $\mathbf{x}_h=\mathcal{J}_h\mathbf{y},$ we have
	$\norm{\mathbf{y}-\mathcal{J}_h\mathbf{y}}_{\delta,\Gamma_C}\leq h^{1-2\delta} \norm{\mathbf{y}}_{3/2-\delta,\Omega}.$ Using the trace estimate (discussed in Section \ref{sec:ModelProblem}) we have $\norm{\bm{\mu}(\mathbf{y})}_{H^{\delta}(\Gamma_C)'}\leq C \big(\norm{\mathbf{y}}_{3/2-\delta,\Omega}+\norm{\bm{\phi}}_{3/2-\delta,\Omega}+\norm{r}_{1/2-\delta,\Omega}\big).$ Putting all these estimates in \eqref{signorini_less_regu} we have
	\begin{align}
		|\rho \,a(\mathbf{y},\mathbf{x}_h-\mathbf{y})+a(\mathbf{y}-\mathbf{x}_h,\bm{\phi})+b(\mathbf{x}_h-\mathbf{y},r)-(\mathbf{u}&-\mathbf{u}_d,\mathbf{y}-\mathbf{x}_h)|\lesssim h^{1-2\delta}.
	\end{align}
	Thus, we have an optimal order(up to the regularity) of convergence of the term \eqref{signori_term}. Hence, Theorems \ref{energy_estm_control}-\ref{estimate of adj pressure} show the optimal order(up to the regularity) of convergence of control, state and adjoint state variables. 	
\end{remark}
\begin{remark}
	For the simplicity of the error analysis we choose the domain $\Omega$ very specific as shown in Figure \ref{example of domain}. But it is clear from the Remark \ref{estm_less_reg} that our error analysis also works for solution with low regularity. This ensures that we can also work with another type of polygonal domain.
\end{remark}

\subsection{Some local $L^2$ and $L^1$ estimate for    \texorpdfstring{$\bm{\mu}(\mathbf{y})$}{\mu(\mathbf{y})} and $\mathbf{y}'$}
\label{Appendix}
\begin{lemma}\label{lemma_req}
	Let $h_e$ be the length of the edge $K\cap \Gamma_{C}$, and $|\mathcal{C}|>0$ and $|\mathcal{N}|>0$. Then the following estimates holds:
	\begin{align}
		\norm{\bm{\mu}(\mathbf{y})}_{0,K\cap \Gamma_C}&\leq\frac{1}{|\mathcal{N}|^{1/2}}\;h_e^{\frac{1}{2}+\delta}\; |\bm{\mu}(\mathbf{y})|_{\delta,K\cap \Gamma_C},\label{lemma_result_1}\\
		\norm{\bm{\mu}(\mathbf{y})}_{L^1(K\cap \Gamma_C)}&\leq\frac{|\mathcal{C}|^{1/2}}{|\mathcal{N}|^{1/2}}\;h_e^{\frac{1}{2}+\delta}\; |\bm{\mu}(\mathbf{y})|_{\delta,K\cap \Gamma_C},\label{lemma_result_2}\\
		\norm{\mathbf{y}'}_{0,K\cap \Gamma_C}&\leq\frac{1}{|\mathcal{C}|^{1/2}}\;h_e^{\frac{1}{2}+\delta}\; |\mathbf{y}'|_{\delta,K\cap \Gamma_C},\label{lemma_result_3}\\
		\norm{\mathbf{y}'}_{L^1(K\cap \Gamma_C)}&\leq\frac{|\mathcal{N}|^{1/2}}{|\mathcal{C}|^{1/2}}\;h_e^{\frac{1}{2}+\delta}\; |\mathbf{y}'|_{\delta,K\cap \Gamma_C},\label{lemma_result_4}
	\end{align}
	where $\bm{\mu}(\mathbf{y})=\rho\frac{\partial \mathbf{y}}{\partial \mathbf{n}}-\frac{\partial \bm{\phi}}{\partial
		\mathbf{n}}-r\mathbf{n}$ and $\mathbf{y}':=(y'_1,y'_2)$ be the tangential derivative of $\mathbf{y}$ on $K\cap\Gamma_C$.
\end{lemma}
\begin{proof}
	Let us start with the $L^2$- estimate of $\bm{\mu}(\mathbf{y})$
	\begin{align*}
		\norm{\bm{\mu}(\mathbf{y})}^2_{0,K\cap \Gamma_C}&=\int_{K\cap \Gamma_C}|\bm{\mu}(\mathbf{y})(\xi)|^2 d\xi\\
		&=\int_{\mathcal{C}}|\bm{\mu}(\mathbf{y})(\xi)|^2 d\xi\quad(\bm{\mu}(\mathbf{y})=0 \;\;\text{on}\;\; \mathcal{N}) \\
		&=\frac{1}{|\mathcal{N}|}\int_{\mathcal{C}}\int_{\mathcal{N}}|\bm{\mu}(\mathbf{y})(\xi)-\bm{\mu}(\mathbf{y})(\eta)|^2d\eta\; d\xi\\
		&\leq\frac{1}{|\mathcal{N}|} \sup_{\mathcal{C}\times \mathcal{N}}|\xi-\eta|^{1+2\delta}\int_{\mathcal{C}}\int_{\mathcal{N}}\frac{|\bm{\mu}(\mathbf{y})(\xi)-\bm{\mu}(\mathbf{y})(\eta)|^2}{|\xi-\eta|^{1+2\delta}}\;d\eta\; d\xi\\
		&\leq\frac{1}{|\mathcal{N}|}\;h_e^{1+2\delta}\; |\bm{\mu}(\mathbf{y})|^2_{\delta,K\cap \Gamma_C},
	\end{align*}
	which proves \eqref{lemma_result_1}. Now we prove \eqref{lemma_result_2} as follows
	\begin{align*}
		\int_{K\cap\Gamma_C}|\bm{\mu}(\mathbf{y})|&=\int_{\mathcal{C}}|\bm{\mu}(\mathbf{y})|,\\
		&\leq |\mathcal{C}|^{1/2} \norm{\bm{\mu}(\mathbf{y})}_{0,\mathcal{C}},\\
		&\leq |\mathcal{C}|^{1/2} \norm{\bm{\mu}(\mathbf{y})}_{0,K\cap\Gamma_C},\\
		&\leq\frac{|\mathcal{C}|^{1/2}}{|\mathcal{N}|^{1/2}}\;h_e^{\frac{1}{2}+\delta}\; |\bm{\mu}(\mathbf{y})|_{\delta,K\cap \Gamma_C}.
	\end{align*}
	The $L^2$- estimate of $\mathbf{y}'$ is derived as follows
	\begin{align*}
		\norm{\mathbf{y}'}^2_{0,K\cap \Gamma_C}&=\int_{K\cap \Gamma_C}|\mathbf{y}'(\xi)|^2 d\xi\\
		&=\int_{\mathcal{N}}|\mathbf{y}'(\xi)|^2 d\xi\quad(\mathbf{y}'=0 \;\;\text{on}\;\; \mathcal{C}) \\
		&=\frac{1}{|\mathcal{C}|}\int_{\mathcal{N}}\int_{\mathcal{C}}|\mathbf{y}'(\xi)-\mathbf{y}'(\eta)|^2d\eta\; d\xi\\
		&\leq\frac{1}{|\mathcal{C}|} \sup_{\mathcal{C}\times \mathcal{N}}|\xi-\eta|^{1+2\delta}\int_{\mathcal{N}}\int_{\mathcal{C}}\frac{|\mathbf{y}'(\xi)-\mathbf{y}'(\eta)|^2}{|\xi-\eta|^{1+2\delta}}\;d\eta\; d\xi\\
		&\leq\frac{1}{|\mathcal{C}|}\;h_e^{1+2\delta}\; |\mathbf{y}'|^2_{\delta,K\cap \Gamma_C},
	\end{align*}
	which proves \eqref{lemma_result_3}. One can easily derive the estimate \eqref{lemma_result_4} by using \eqref{lemma_result_3}.
\end{proof}

\section{Numerical Experiments}\label{sec:Numerics}
In this section, we are going to validate the a priori error
estimates for the error in the control, state, and adjoint state numerically. Here we consider two model examples with known exact solutions. For the numerical experiments we slightly modify the optimal control problem which is as follows:
$$\text{minimize}\;  \tilde{J}(\mathbf{w},\mathbf{x})=\frac{1}{2}\|\mathbf{w}-\mathbf{u}_d\|_{0,\Omega}^2+\frac{\rho}{2}\|\nabla (\mathbf{x}-\mathbf{y}_d)\|_{0,\Omega}^2,$$
subject to the PDE,
\begin{subequations}\label{stokes_eq_intro_1}
	\begin{align}
	-\Delta \mathbf{w}+\nabla{p}&=\mathbf{f} \quad \text{in}\;\Omega, \\
	\nabla\cdot{\mathbf{w}}&=0 \quad \text{in}\; \Omega,\nonumber\\
	\bf{w}&=\mathbf{x} \quad \text{on}\; \Gamma_C,\nonumber\\
	\bf{w}&=\mathbf{0} \quad \text{on}\; \Gamma_D,\nonumber
	\end{align}
\end{subequations}
the control set is given by
$$\mathbf{Q}_{ad}:=\{\mathbf{x}\in \mathbf{H}^1(\Omega): \bm{\gamma}_{0}(\mathbf{x})=\mathbf{0}\text{ on } \Gamma_D, \mathbf{y}_a\leq\bm{\gamma}_{0}(\mathbf{x})\leq \mathbf{y}_b\text{ a.e. on } \Gamma_C\},$$
where the function $\mathbf{y}_d$ is given and the boundary $\partial\Omega=\Gamma_C\cup \bar{\Gamma}_D$. Consequently, the discrete optimality system takes the form
 \begin{subequations}
 	\begin{align*}
 		\mathbf{u}_h&=\mathbf{w}_h+\mathbf{y}_h,\quad \mathbf{w}_h\in \mathbf{V}_h,\\
 		a(\mathbf{w}_h,\mathbf{z}_h)+b(\mathbf{z}_h,p_H) &=( \mathbf{f},\mathbf{z}_h)-a(\mathbf{y}_h,\mathbf{z}_h) \;\;\;
 		{\rm for~all}\;\mathbf{z}_h \in \mathbf{V}_h,\\ 
 		b(\mathbf{u}_h,q_H)&=0\;\quad {\rm for~all}\; q_H \in M_H ,\\
 		a(\mathbf{z}_h,\bm{\phi}_h)-b(\mathbf{z}_h,r_H) &=( \mathbf{u}_h-\mathbf{u_d},\mathbf{z}_h) \;\;\;
 		{\rm for~all}\;\mathbf{z}_h \in \mathbf{V}_h ,\\ 
 		b(\bm{\phi}_h,q_H)&=0\;\quad {\rm for~all}\; q_H \in M_H ,\\
 		\rho\, a(\mathbf{y}_h,\mathbf{x}_h-\mathbf{y}_h)\geq& a(\mathbf{x}_h-\mathbf{y}_h,\bm{\phi}_h)-b(\mathbf{x}_h-\mathbf{y}_h,r_H)\\&-(\mathbf{u}_h-\mathbf{u_d},\mathbf{x}_h-\mathbf{y}_h)\;\quad {\rm for~all}\; \mathbf{x}_h\in
 		\mathbf{Q}^h_{ad},
 	\end{align*}
 \end{subequations}
where,
 $\mathbf{V}_h:=\{\mathbf{v}_h\in \mathbf{H}^1_{0}(\Omega): \mathbf{v}_h|_T\in \mathbf{P}_1(K)\;\; \forall K \in \mathcal{T}_h\}$
 and $\mathbf{Q}^h_{ad}=\mathbf{Q}_h\cap\mathbf{Q}_{ad}.$ The set $\mathbf{Q}_h$ is defined by
 $\mathbf{Q}_h=\{\mathbf{x}_h \in \mathbf{H}^1(\Omega): \bm{\gamma}_{0}(\mathbf{x})=\mathbf{0}\text{ on } \Gamma_D,\; \mathbf{x}_h|_T\in \mathbf{P}_1(K),\;\; \forall K \in \mathcal{T}_h\}.$

\begin{figure}
 	\centering
 	\begin{subfigure}[t]{0.4\textwidth}
 		\centering
 		\includegraphics[width=8cm,height=2.5in]{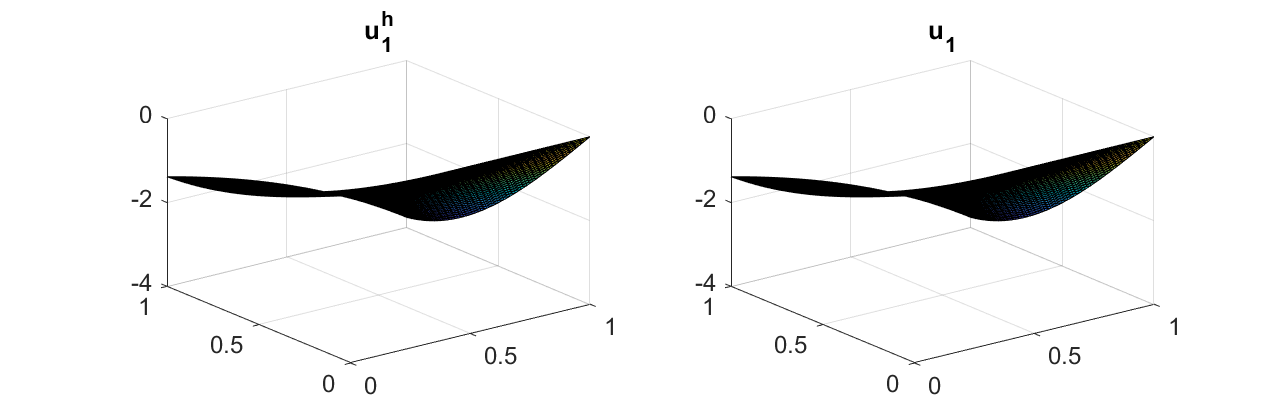}
 		\caption{}
 	\end{subfigure}
 	\hfill
 	\begin{subfigure}[t]{0.4\textwidth}
 		\centering
 		\includegraphics[width=7cm,height=2.5in]{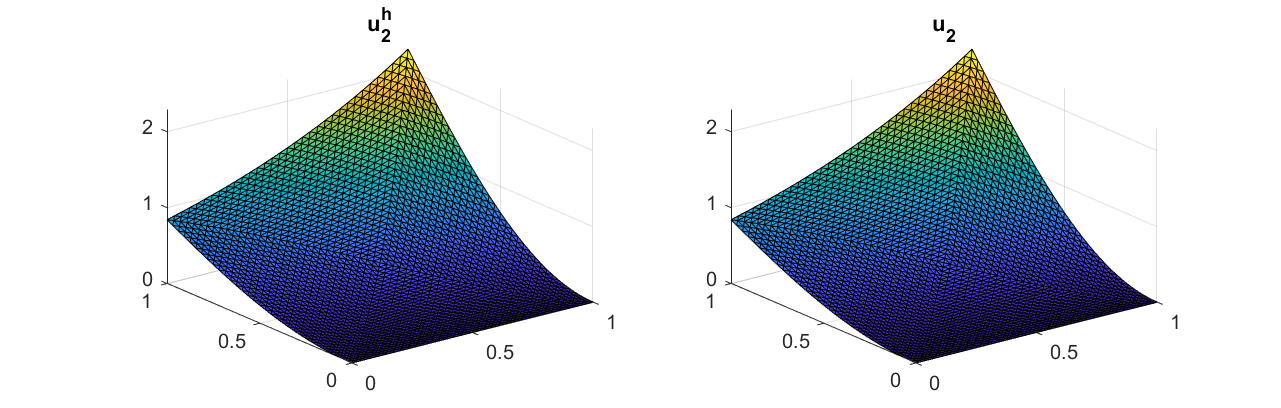}
 		\caption{}
 	\end{subfigure}
 		\caption{(A) Discrete state $(\mathbf{u}_1^h)$ and continuous state $(\mathbf{u}_1)$ for Example \ref{Ex.1.1} (B) Discrete state $(\mathbf{u}_2^h)$ and continuous state $(\mathbf{u}_2)$  for Example \ref{Ex.1.1}}
 		\label{fig:state}
 \end{figure}

\begin{figure}
	\centering
	\begin{subfigure}[b]{0.4\textwidth}
		\centering
		\includegraphics[width=8cm,height=2.5in]{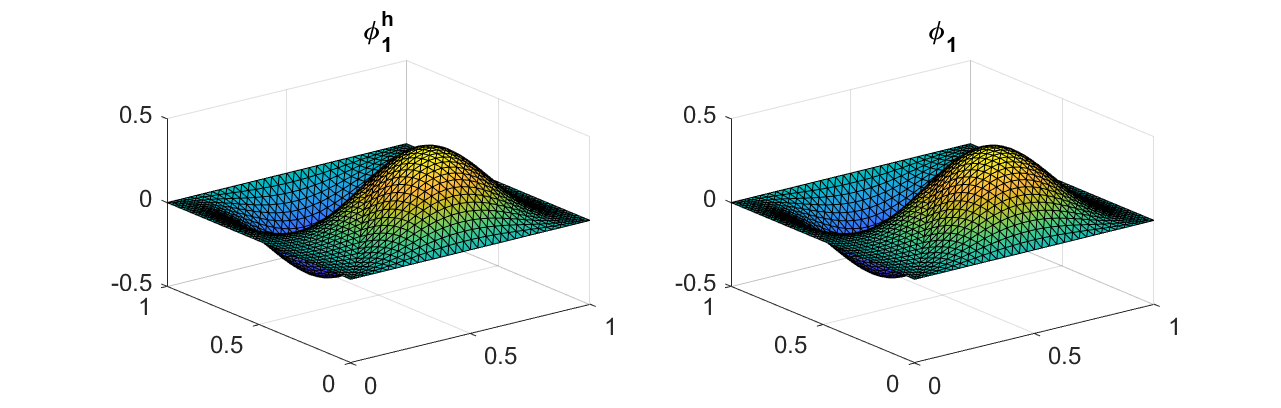}
		\caption{}
	\end{subfigure}
	\hfill
	\begin{subfigure}[b]{0.4\textwidth}
		\centering
		\includegraphics[width=7cm,height=2.5in]{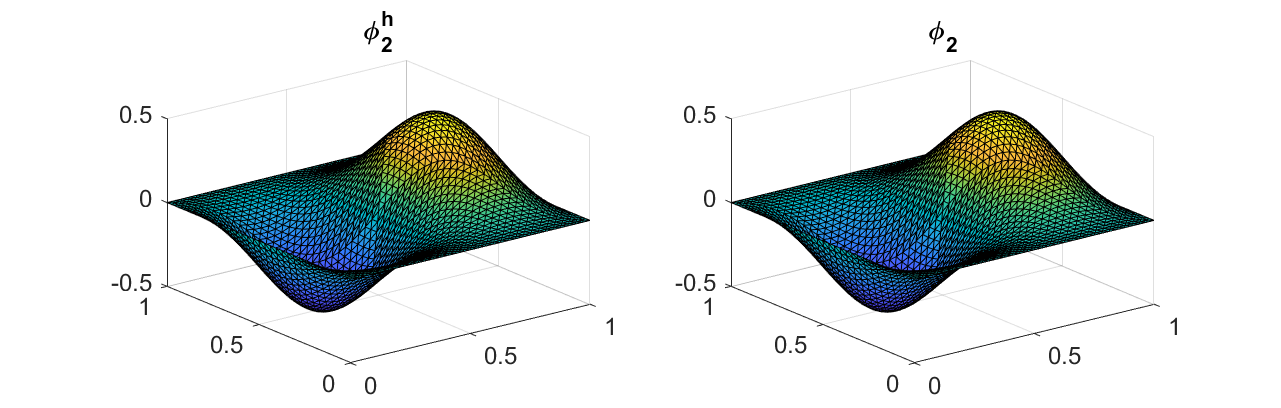}
		\caption{}
	\end{subfigure}
		\caption{(A) Discrete adj-state $(\bm{\phi}_1^h)$ and continuous adj-state $(\bm{\phi}_1)$ for Example \ref{Ex.1.1} (B) Discrete adj-state $(\bm{\phi}_2^h)$ and continuous adj-state $(\bm{\phi}_2)$ for Example \ref{Ex.1.1}}
		\label{fig:adj_state}
\end{figure}

\begin{figure}
	\centering
	\begin{subfigure}[b]{0.4\textwidth}
		\centering
		\includegraphics[width=8cm,height=2.5in]{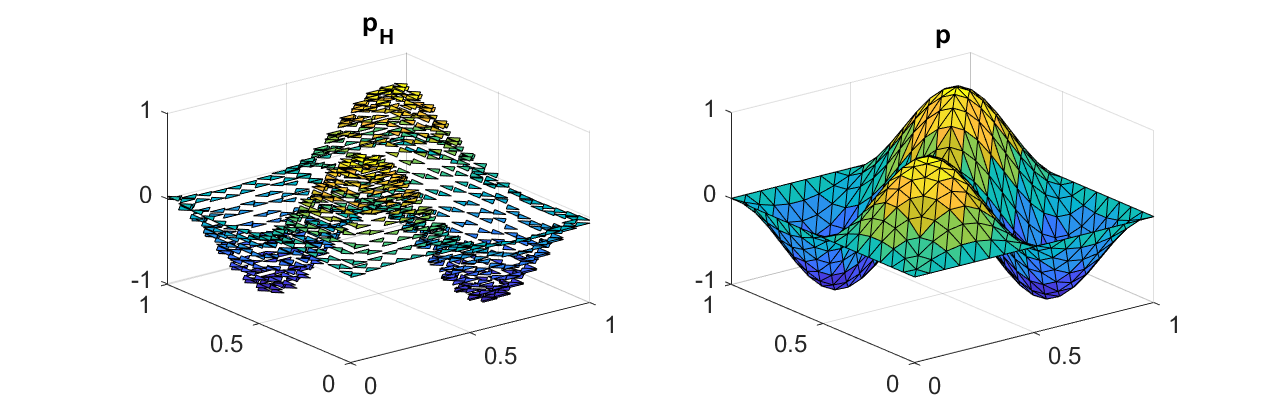}
		\caption{}
	\end{subfigure}
	\hfill
	\begin{subfigure}[b]{0.4\textwidth}
		\centering
		\includegraphics[width=7cm,height=2.5in]{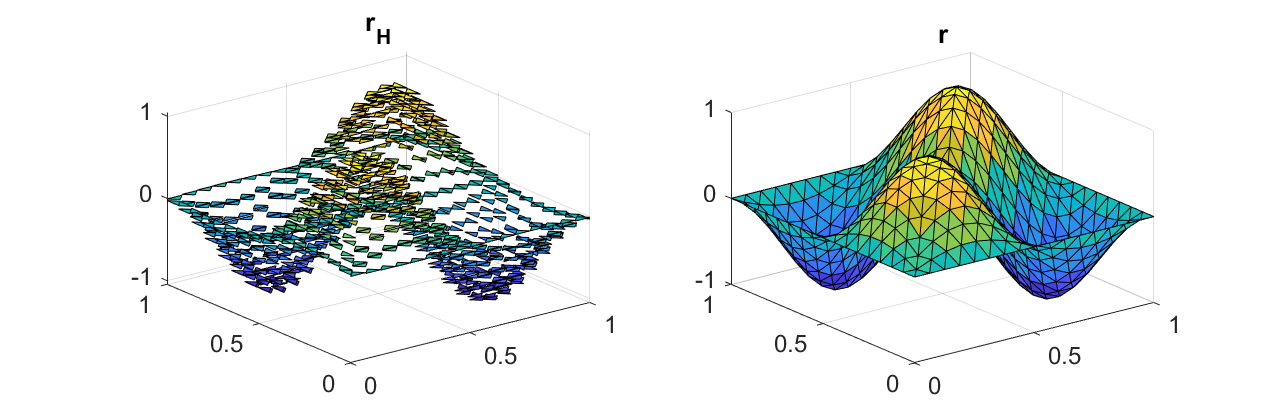}
		\caption{}
	\end{subfigure}
		\caption{(A) Discrete pressure $(p_H)$ and continuous pressure $(p)$ for Example \ref{Ex.1.1} (B) Discrete adj-pressure $(r_H)$ and continuous adj-pressure $(r)$  for Example \ref{Ex.1.1} }
		\label{fig:pressure}
\end{figure}

\begin{figure}
 	\centering
 	\begin{subfigure}[b]{0.4\textwidth}
 		\centering
 		\includegraphics[width=8cm,height=2.5in]{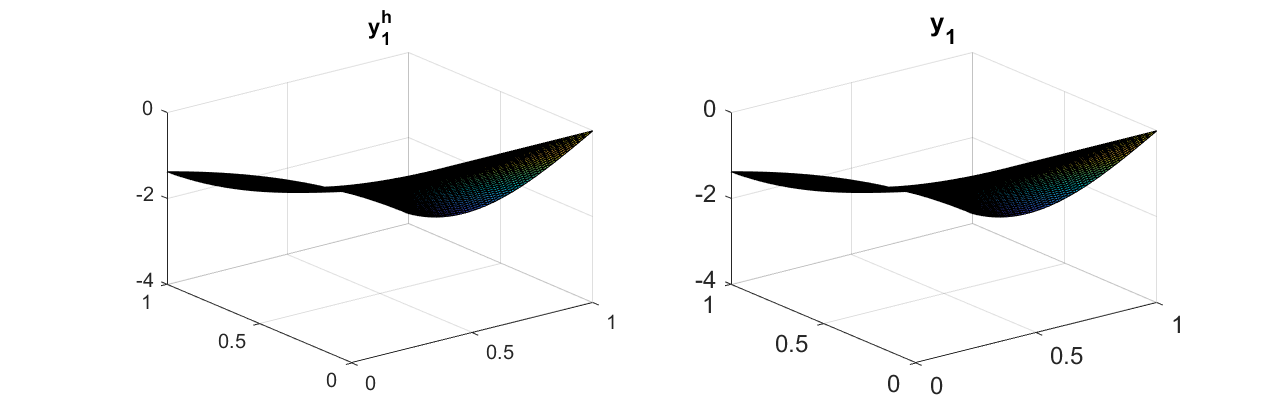}
 		\caption{}
 	\end{subfigure}
 	\hfill
 	\begin{subfigure}[b]{0.4\textwidth}
 		\centering
 		\includegraphics[width=7cm,height=2.5in]{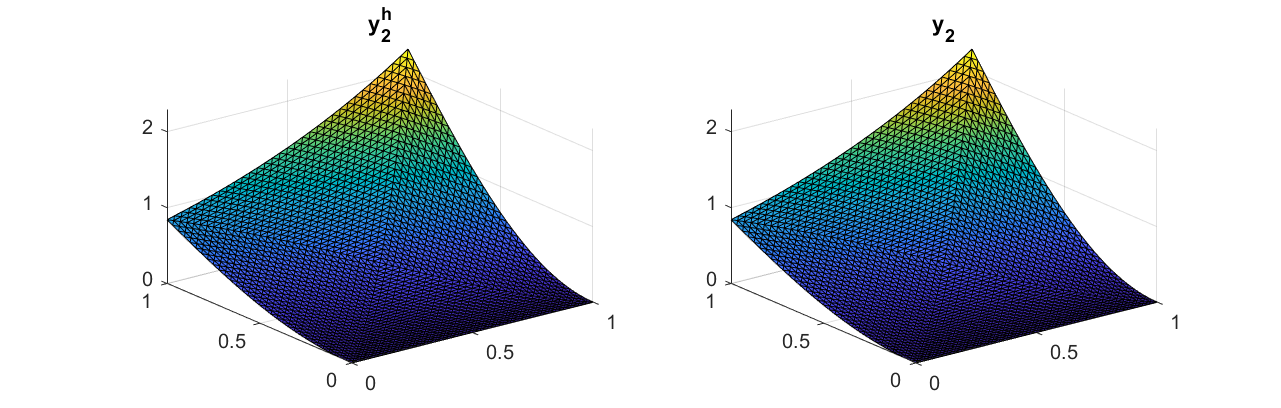}
 		\caption{}
 	\end{subfigure}
 		\caption{(A) Discrete control $(\mathbf{y}_1^h)$ and continuous control $(\mathbf{y}_1)$ for Example \ref{Ex.1.1} (B)Discrete control $(\mathbf{y}_2^h)$ and continuous control $(\mathbf{y}_2)$  for Example \ref{Ex.1.1} }
 		\label{fig:control}
 \end{figure}

 \begin{example}\label{Ex.1.1}
 	Let the computational domain $\Omega=(0,1)^2$, and $\Gamma_D=(0,1)\times \{0\},$ $\Gamma_C=\partial \Omega \backslash \Gamma_D$. We choose the constants $\rho=10^{-2},$ $\mathbf{y}_a=(-4,0),$ and $\mathbf{y}_b=(0,2.5).$ The state and adjoint state variables are given by
 	\begin{align}
 	{\bf u}={\bf y}= \left(
 	\begin{array}{c}
 	-\exp(x)(y\cos(y)+\sin(y))\\
 \exp(x)y\sin(y)
 	\end{array}
 	\right),~~~ p = \sin(2\pi x)\sin(2 \pi y),
 	\end{align}
 	and
 	\begin{align}
 	\bm{\phi}= \left(
 	\begin{array}{c}
 	(\sin(\pi x))^2 \sin(\pi y)\cos(\pi y)\\
 	-(\sin(\pi y))^2 \sin(\pi x) \cos(\pi x)
 	\end{array}
 	\right),~~~ r = \sin(2\pi x)\sin(2 \pi y).
 	\end{align}
 	We choose $\mathbf{u}$ and $\bm{\phi}$ such that $\nabla\cdot{\mathbf{u}}=\nabla\cdot{\bm{\phi}}=0 \quad \text{in}\; \Omega$  and $\bm{\phi}=\mathbf{0}~~ \text{on}~~ \partial\Omega.$ The data of the problem are  chosen such that $\mathbf{f} = -\Delta \mathbf{u}+\nabla{p},~\mathbf{u}_d= \mathbf{u}+\Delta \bm{\phi}+\nabla{r}~ \text{and}~\mathbf{y}_d=\mathbf{y}$.

 	
 \end{example}
The discrete solution is computed on several uniform grids with mesh sizes $h=\frac{1}{2^i}, i=2,...,6$ for the velocity variable and mesh sizes $H=2h$ for the pressure variable. To solve the optimal control problem numerically, we have used the primal-dual active set algorithm. The continuous and discrete approximations for the state velocity variables using conforming $\mathbf{P}_1$(in fine mesh) finite elements are shown in Figure \ref{fig:state}. The continuous and
 discrete approximations for the pressure variables using conforming $P_0$(in coarse mesh) finite elements are shown in Figure \ref{fig:pressure}(A). The continuous and
 discrete approximations for the  adjoint state velocity variables using conforming $\mathbf{P}_1$(in fine mesh) finite elements are shown in Figure \ref{fig:adj_state}. The continuous and
 discrete approximations for the pressure variables using conforming $P_0$(in coarse mesh) finite elements are shown in Figure \ref{fig:pressure}(B). The continuous and
 discrete approximations for the control variables using conforming $\mathbf{P}_1$(in fine mesh) finite elements are shown in Figure \ref{fig:control}.

 Tables \ref{table1.1} and \ref{table1.2} show the computed errors and orders of convergence of the state and adjoint state variables respectively for the Example \ref{Ex.1.1}. The errors and orders of convergence  of the control variable is shown in Table \ref{table1.3}. The numerical convergence rates with respect to the energy norm for the state, adjoint state and control variables are linear as predicted theoretically.

  \begin{table}[h!]
  	\begin{center}
  		\footnotesize
  		\vspace{0.2cm}
  		\caption{Energy errors and convergence rates of the state variable for Example \ref{Ex.1.1}.}
  		\label{table1.1}
  		\begin{tabular}{ |c|c|c|c|c|c|}
  			\hline
  			$h$ & $\norm{\nabla(\bf{u}-\bf{u}_h)}_{0,\Omega}$ & order & $H$ &  $\norm{p-p_{H}}_{0,\Omega}$ & order\\
  			
  			\hline
  			0.2500  &  1.2670   &     0    &  0.5000  &  0.3327   &      0\\
  			0.1250  & 0.7124    &  0.8307  &  0.2500  &  0.7025   &-1.0782\\
  			0.0625  &   0.3502  & 1.0243   &   0.1250 &   0.3584 &   0.9709\\
  			0.0312  &  0.1770   & 0.9843   &  0.0625  &  0.1878  &  0.9323\\
  			0.0156  & 0.0888    &  0.9955  & 0.0312   & 0.0949   & 0.9852\\

  			\hline
  		\end{tabular}
  	\end{center}
  \end{table}

  \begin{table}[h!]
  	\begin{center}
  		\footnotesize
  		\vspace{0.2cm}
  		\caption{Energy errors and convergence rates of the adjoint state variable for Example \ref{Ex.1.1}.}
  		\label{table1.2}
  		\begin{tabular}{ |c|c|c|c|c|c|}
  			\hline
  			$h$ & $\norm{\nabla(\bm{\phi}-\bm{\phi}_h)}_{0,\Omega}$ & order & $H$ &  $\norm{r-r_{H}}_{0,\Omega}$ & order\\
  			
  			\hline
  			
  			0.2500 &   2.3420 &        0&    0.5000 &   0.0153&         0\\
  			0.1250  &  1.3102  &  0.8380 &   0.2500  &  0.7643 &  -5.6448\\
  			0.0625   & 0.6934   & 0.9181  &  0.1250   & 0.3786  &  1.0136\\
  			0.0312    &0.3537    &0.9711   & 0.0625    &0.1975   & 0.9386\\
  			0.0156    &0.1777    &0.9928    &0.0312    &0.0995    &0.9890\\

  			\hline
  		\end{tabular}
  	\end{center}
  \end{table}

   \begin{table}[h!]
   	\begin{center}
   		\footnotesize
   		\vspace{0.2cm}
   		\caption{Energy errors and convergence rates of the control variable for Example \ref{Ex.1.1}.}
   		\label{table1.3}
   		\begin{tabular}{ |c|c|c|}
   			\hline
   			$h$ & $\norm{\nabla(\bf{y}-\bf{y}_h)}_{0,\Omega}$ & order\\
   			
   			\hline
   			0.2500&    1.2307&         0\\
   			0.1250 &   0.6167 &   0.9969\\
   			0.0625  &  0.3085  &  0.9993\\
   			0.0312   & 0.1543   & 0.9998\\
   			0.0156    &0.0771    &1.0000\\

   			\hline
   		\end{tabular}
   	\end{center}
   \end{table}

  \begin{example}\label{Ex.1.2}
  	In this example, we report the results of numerical tests
  	carried out for the $L$-shaped domain $\Omega=(-\frac{1}{2},\frac{1}{2})^2\backslash\big((0,\frac{1}{2}) \times (0,-\frac{1}{2})\big),$ with the boundaries $\Gamma_D=(0,\frac{1}{2})\times \{0\}\cup \{0\} \times (0,-\frac{1}{2})$ and $\Gamma_C=\partial \Omega \backslash \Gamma_D$. We choose the constants $\rho=10^{-2},$ $\mathbf{y}_a=(-0.6,-0.6)$ and $\mathbf{y}_b=(0.6,0.6)$. The state and adjoint state variables are given by
  	\begin{align}
  	{\bf u}={\bf y}= \left(
  	\begin{array}{c}
  	(\sin(\pi x))^2 \sin(\pi y)\cos(\pi y)\\
  		-(\sin(\pi y))^2 \sin(\pi x) \cos(\pi x)
  	\end{array}
  	\right),~~~ p = \sin(2\pi x)\sin(2 \pi y),
  	\end{align}
  and
  	\begin{align}
  	\bm{\phi}= \left(
  	\begin{array}{c}
  	(\sin(2\pi x))^2 \sin(2\pi y)\cos(2\pi y)\\
  	-(\sin(2\pi y))^2 \sin(2\pi x) \cos(2\pi x)
  	\end{array}
  	\right),~~~ r = \sin(2\pi x)\sin(2 \pi y).
  	\end{align}
  We choose $\mathbf{u}$ and $\bm{\phi}$ such that $\nabla\cdot{\mathbf{u}}=\nabla\cdot{\bm{\phi}}=0 \quad \text{in}\; \Omega$ and $\bm{\phi}=0~~ \text{on}~~ \partial\Omega.$
  	The data of the problem are chosen such that $\mathbf{f} = -\Delta \mathbf{u}+\nabla{p},~\mathbf{u}_d= \mathbf{u}+\Delta \bm{\phi}+\nabla{r},~{\rm and}~ \mathbf{y}_d=\mathbf{y}$.

  	
  \end{example}
In the same way as in Example \ref{Ex.1.1}, we produce a sequence of meshes. Tables \ref{table1.4} and \ref{table1.5}  show the computed errors and orders of convergence  of the state and adjoint state variables respectively for Example \ref{Ex.1.2}. The errors and orders of convergence  of the control variable is shown in Table \ref{table1.6}. The observed orders of convergences agree with the predicted ones as seen in the tables.

  \begin{table}[h!]
  	\begin{center}
  		\footnotesize
  		\vspace{0.2cm}
  		\caption{Energy errors and convergence rates of the state variable for Example \ref{Ex.1.2}.}
  		\label{table1.4}
  		\begin{tabular}{ |c|c|c|c|c|c|}
  			\hline
  			$h$ & $\norm{\nabla(\bf{u}-\bf{u}_h)}_{0,\Omega}$ & order & $H$ &  $\norm{p-p_{H}}_{0,\Omega}$ & order\\
  			
  			\hline
  			0.1768&    1.6516&         0 &   0.3536&    0.2711&         0\\
  			0.0884 &   0.9670 &   0.7723  &  0.1768 &   0.5179 &   -0.9338\\
  			0.0442  &  0.5029  &  0.9433   & 0.0884  &  0.2799  &  0.8875\\
  			0.0221   & 0.2539   & 0.9862    &0.0442   & 0.1473   & 0.9259\\
  			0.0110    &0.1272    &0.9966    &0.0221    &0.0771    &0.9335\\

  			\hline
  		\end{tabular}
  	\end{center}
  \end{table}

  \begin{table}[h!]
  	\begin{center}
  		\footnotesize
  		\vspace{0.2cm}
  		\caption{Energy errors and convergence rates of the adjoint state variable for the Example \ref{Ex.1.2}.}
  		\label{table1.5}
  		\begin{tabular}{ |c|c|c|c|c|c|}
  			\hline
  			$h$ & $\norm{\nabla(\bm{\phi}-\bm{\phi}_h)}_{0,\Omega}$ & order & $H$ &  $\norm{r-r_{H}}_{0,\Omega}$ & order\\
  			
  			\hline
  			0.1768&    5.4417&         0&    0.3536&    0.0844&        0\\
  			0.0884 &   3.0884 &   0.8172 &   0.1768 &   0.1857  & -1.1376\\
  			0.0442  &  1.8577&    0.7333  &  0.0884  &  0.0987   &0.9107\\
  			0.0221   & 0.9813 &   0.9208   & 0.0442   & 0.0511    &0.9486\\
  			0.0110    &0.4979  &  0.9788    &0.0221    &0.0262    &0.9634\\

  			\hline
  		\end{tabular}
  	\end{center}
  \end{table}

  \begin{table}[h!]
  	\begin{center}
  		\footnotesize
  		\vspace{0.2cm}
  		\caption{Energy errors and convergence rates of the control variable for the Example \ref{Ex.1.2}.}
  		\label{table1.6}
  		\begin{tabular}{ |c|c|c|}
  			\hline
  			$h$ & $\norm{\nabla(\bf{y}-\bf{y}_h)}_{0,\Omega}$ & order\\
  			
  			\hline
  			
  			0.1768&    1.6656&         0\\
  			0.0884 &   0.9364 &   0.8308\\
  			0.0442  &  0.4827  &  0.9560\\
  			0.0221   & 0.2432   & 0.9889\\
  			0.0110    &0.1218    &0.9972\\

  			\hline
  		\end{tabular}
  	\end{center}
  \end{table}

\section*{Acknowledgment}
The authors would like to sincerely thank Professor J.-P. Raymond for some fruitful discussions on the regularity of the Stokes solution.

\bibliography{bibliography}

\bibliographystyle{abbrvnat}
\end{document}